
\magnification 1200
\hsize 13.2cm
\vsize 20cm
\parskip 3pt plus 1pt
\parindent 5mm

\def\\{\hfil\break}


\font\seventeenbf=cmbx10 at 17.28pt

\font\twelvebf=cmbx10 at 12pt
\font\eightbf=cmbx8
\font\sixbf=cmbx6

\font\eighti=cmmi8
\font\sixi=cmmi6

\font\eightrm=cmr8
\font\sixrm=cmr6

\font\eightsy=cmsy8
\font\sixsy=cmsy6

\font\eightit=cmti8
\font\eighttt=cmtt8
\font\eightsl=cmsl8

\font\seventeenbsy=cmbsy10 at 17.28pt

\font\twelvebsy=cmbsy10 at 12pt
\font\tenbsy=cmbsy10
\font\eightbsy=cmbsy8
\font\sevenbsy=cmbsy7
\font\sixbsy=cmbsy6
\font\fivebsy=cmbsy5

\font\tenmsa=msam10

\font\sevenmsa=msam7
\font\fivemsa=msam5
\newfam\msafam
  \textfont\msafam=\tenmsa
  \scriptfont\msafam=\sevenmsa
  \scriptscriptfont\msafam=\fivemsa

\font\tenmsb=msbm10
\font\eightmsb=msbm8
\font\sevenmsb=msbm7
\font\fivemsb=msbm5
\newfam\msbfam
  \textfont\msbfam=\tenmsb
  \scriptfont\msbfam=\sevenmsb
  \scriptscriptfont\msbfam=\fivemsb
\def\Bbb{\fam\msbfam\tenmsb}

\font\tenCal=eusm10
\font\sevenCal=eusm7
\font\fiveCal=eusm5
\newfam\Calfam
  \textfont\Calfam=\tenCal
  \scriptfont\Calfam=\sevenCal
  \scriptscriptfont\Calfam=\fiveCal
\def\Cal{\fam\Calfam\tenCal}

\font\teneuf=eusm10
\font\teneuf=eufm10
\font\seveneuf=eufm7
\font\fiveeuf=eufm5
\newfam\euffam
  \textfont\euffam=\teneuf
  \scriptfont\euffam=\seveneuf
  \scriptscriptfont\euffam=\fiveeuf

\font\seventeenbfit=cmmib10 at 17.28pt

\font\twelvebfit=cmmib10 at 12pt
\font\tenbfit=cmmib10
\font\eightbfit=cmmib8
\font\sevenbfit=cmmib7
\font\sixbfit=cmmib6
\font\fivebfit=cmmib5
\newfam\bfitfam
  \textfont\bfitfam=\tenbfit
  \scriptfont\bfitfam=\sevenbfit
  \scriptscriptfont\bfitfam=\fivebfit


\catcode`\@=11
\def\eightpoint{%
  \textfont0=\eightrm \scriptfont0=\sixrm \scriptscriptfont0=\fiverm
  \def\rm{\fam\z@\eightrm}%
  \textfont1=\eighti \scriptfont1=\sixi \scriptscriptfont1=\fivei
  \def\oldstyle{\fam\@ne\eighti}%
  \textfont2=\eightsy \scriptfont2=\sixsy \scriptscriptfont2=\fivesy
  \textfont\itfam=\eightit
  \def\it{\fam\itfam\eightit}%
  \textfont\slfam=\eightsl
  \def\sl{\fam\slfam\eightsl}%
  \textfont\bffam=\eightbf \scriptfont\bffam=\sixbf
  \scriptscriptfont\bffam=\fivebf
  \def\bf{\fam\bffam\eightbf}%
  \textfont\ttfam=\eighttt
  \def\tt{\fam\ttfam\eighttt}%
  \textfont\msbfam=\eightmsb
  \def\Bbb{\fam\msbfam\eightmsb}%
  \abovedisplayskip=9pt plus 2pt minus 6pt
  \abovedisplayshortskip=0pt plus 2pt
  \belowdisplayskip=9pt plus 2pt minus 6pt
  \belowdisplayshortskip=5pt plus 2pt minus 3pt
  \smallskipamount=2pt plus 1pt minus 1pt
  \medskipamount=4pt plus 2pt minus 1pt
  \bigskipamount=9pt plus 3pt minus 3pt
  \normalbaselineskip=9pt
  \setbox\strutbox=\hbox{\vrule height7pt depth2pt width0pt}%
  \let\bigf@ntpc=\eightrm \let\smallf@ntpc=\sixrm
  \normalbaselines\rm}
\catcode`\@=12

\def\eightpointbf{%
 \textfont0=\eightbf   \scriptfont0=\sixbf   \scriptscriptfont0=\fivebf
 \textfont1=\eightbfit \scriptfont1=\sixbfit \scriptscriptfont1=\fivebfit
 \textfont2=\eightbsy  \scriptfont2=\sixbsy  \scriptscriptfont2=\fivebsy
 \eightbf
 \baselineskip=10pt}

\def\tenpointbf{%
 \textfont0=\tenbf   \scriptfont0=\sevenbf   \scriptscriptfont0=\fivebf
 \textfont1=\tenbfit \scriptfont1=\sevenbfit \scriptscriptfont1=\fivebfit
 \textfont2=\tenbsy  \scriptfont2=\sevenbsy  \scriptscriptfont2=\fivebsy
 \tenbf}

\def\twelvepointbf{%
 \textfont0=\twelvebf   \scriptfont0=\eightbf   \scriptscriptfont0=\sixbf
 \textfont1=\twelvebfit \scriptfont1=\eightbfit \scriptscriptfont1=\sixbfit
 \textfont2=\twelvebsy  \scriptfont2=\eightbsy  \scriptscriptfont2=\sixbsy
 \twelvebf
 \baselineskip=14.4pt}

\def\seventeenpointbf{%
 \textfont0=\seventeenbf  \scriptfont0=\twelvebf  \scriptscriptfont0=\eightbf
 \textfont1=\seventeenbfit\scriptfont1=\twelvebfit\scriptscriptfont1=\eightbfit
 \textfont2=\seventeenbsy \scriptfont2=\twelvebsy \scriptscriptfont2=\eightbsy
 \seventeenbf
 \baselineskip=20.736pt}


\newdimen\srdim \srdim=\hsize
\newdimen\irdim \irdim=\hsize
\def\NOSECTREF#1{\noindent\hbox to \srdim{\null\dotfill ???(#1)}}
\def\SECTREF#1{\noindent\hbox to \srdim{\csname REF\romannumeral#1\endcsname}}
\def\INDREF#1{\noindent\hbox to \irdim{\csname IND\romannumeral#1\endcsname}}
\newlinechar=`\^^J
\def\openauxfile{
  \immediate\openin1\jobname.aux
  \ifeof1
  \message{^^JCAUTION\string: you MUST run TeX a second time^^J}
  \let\sectref=\NOSECTREF \let\indref=\NOSECTREF
  \else
  \input \jobname.aux
  \message{^^JCAUTION\string: if the file has just been modified you may
    have to run TeX twice^^J}
  \let\sectref=\SECTREF \let\indref=\INDREF
  \fi
  \message{to get correct page numbers displayed in Contents or Index
    Tables^^J}
  \immediate\openout1=\jobname.aux
  \let\END=\end \def\end{\immediate\closeout1\END}}

\newbox\titlebox   \setbox\titlebox\hbox{\hfil}
\newbox\sectionbox \setbox\sectionbox\hbox{\hfil}
\def\folio{\ifnum\pageno=1 \hfil \else \ifodd\pageno
           \hfil {\eightpoint\copy\sectionbox\kern8mm\number\pageno}\else
           {\eightpoint\number\pageno\kern8mm\copy\titlebox}\hfil \fi\fi}
\footline={\hfil}
\headline={\folio}

\def\titlerunning#1{\setbox\titlebox\hbox{\eightpoint #1}}
\def\title#1{\noindent\hfil$\smash{\hbox{\seventeenpointbf #1}}$\hfil
             \titlerunning{#1}\medskip}

\newcount\numbersection \numbersection=-1
\def\sectionrunning#1{\setbox\sectionbox\hbox{\eightpoint #1}
  \immediate\write1{\string\def \string\REF
      \romannumeral\numbersection \string{%
      \noexpand#1 \string\dotfill \space \number\pageno \string}}}
\def\section#1{%
  \par\vskip0.666cm\penalty -100
  \vbox{\baselineskip=14.4pt\noindent{{\twelvepointbf #1}}}
  \vskip2pt
  \penalty 500
  \advance\numbersection by 1
  \sectionrunning{#1}}

\def\subsection#1{%
  \par\vskip0.5cm\penalty -100
  \vbox{\noindent{{\tenpointbf #1}}}
  \vskip1pt
  \penalty 500}

\newcount\numberindex \numberindex=0
\def\index#1#2{%
  \advance\numberindex by 1
  \immediate\write1{\string\def \string\IND #1%
     \romannumeral\numberindex \string{%
     \noexpand#2 \string\dotfill \space \string\S \number\numbersection,
     p.\string\ \space\number\pageno \string}}}

\newdimen\itemindent \itemindent=\parindent

\def\item#1{\par\noindent\hangindent\itemindent%
            \rlap{#1}\kern\itemindent\ignorespaces}
\def\itemitem#1{\par\noindent\hangindent2\itemindent%
            \kern\itemindent\rlap{#1}\kern\itemindent\ignorespaces}
\def\itemitemitem#1{\par\noindent\hangindent3\itemindent%
            \kern2\itemindent\rlap{#1}\kern\itemindent\ignorespaces}

\long\def\claim#1|#2\endclaim{\par\vskip 5pt\noindent
{\tenpointbf #1.}\ {\it #2}\par\vskip 5pt}

\def\proof{\noindent{\it Proof}}
\def\hotimes{\hat \otimes}

\def\today{\ifcase\month\or
January\or February\or March\or April\or May\or June\or July\or August\or
September\or October\or November\or December\fi \space\number\day,
\number\year}

\catcode`\@=11
\newcount\@tempcnta \newcount\@tempcntb
\def\timeofday{{%
\@tempcnta=\time \divide\@tempcnta by 60 \@tempcntb=\@tempcnta
\multiply\@tempcntb by -60 \advance\@tempcntb by \time
\ifnum\@tempcntb > 9 \number\@tempcnta:\number\@tempcntb
  \else\number\@tempcnta:0\number\@tempcntb\fi}}
\catcode`\@=12

\def\bibitem#1&#2&#3&#4&%
{\hangindent=1.8cm\hangafter=1
\noindent\rlap{\hbox{\eightpointbf #1}}\kern1.8cm{\rm #2}{\it #3}{\rm #4.}}


\def\bC{{\Bbb C}}

\def\bN{{\Bbb N}}

\def\bQ{{\Bbb Q}}
\def\bR{{\Bbb R}}

\def\bZ{{\Bbb Z}}

\def\cC{{\Cal C}}
\def\cE{{\Cal E}}
\def\cF{{\Cal F}}

\def\cH{{\cal H}}
\def\cI{{\Cal I}}

\def\cO{{\Cal O}}
\def\cR{{\Cal R}}
\def\cS{{\Cal S}}
\def\cA{{\Cal A}}
\def\cB{{\Cal B}}
\def\cS{{\Cal S}}
\def\cT{{\Cal T}}
\def\cX{{\Cal X}}


\def\Rangle{\rangle\kern-2pt\rangle}
\def\Langle{\langle\kern-2pt\langle}

\def\square{{\hfill \hbox{
\vrule height 1.453ex  width 0.093ex  depth 0ex
\vrule height 1.5ex  width 1.3ex  depth -1.407ex\kern-0.1ex
\vrule height 1.453ex  width 0.093ex  depth 0ex\kern-1.35ex
\vrule height 0.093ex  width 1.3ex  depth 0ex}}}
\def\qed{\kern10pt$\square$}
\def\hexnbr#1{\ifnum#1<10 \number#1\else
 \ifnum#1=10 A\else\ifnum#1=11 B\else\ifnum#1=12 C\else
 \ifnum#1=13 D\else\ifnum#1=14 E\else\ifnum#1=15 F\fi\fi\fi\fi\fi\fi\fi}
\def\msatype{\hexnbr\msafam}
\def\msbtype{\hexnbr\msbfam}
\mathchardef\restriction="3\msatype16   
\mathchardef\boxsquare="3\msatype03
\mathchardef\preccurlyeq="3\msatype34
\mathchardef\compact="3\msatype62
\mathchardef\smallsetminus="2\msbtype72   \let\ssm\smallsetminus
\mathchardef\subsetneq="3\msbtype28
\mathchardef\supsetneq="3\msbtype29
\mathchardef\leqslant="3\msatype36   \let\le\leqslant
\mathchardef\geqslant="3\msatype3E   \let\ge\geqslant
\mathchardef\stimes="2\msatype02
\mathchardef\ltimes="2\msbtype6E
\mathchardef\rtimes="2\msbtype6F

\def\dbar{\overline\partial}
\def\ddbar{\partial\overline\partial}


\let\la=\longrightarrow

\let\text=\hbox
\def\buildo#1^#2{\mathop{#1}\limits^{#2}}
\def\buildu#1_#2{\mathop{#1}\limits_{#2}}
\def\ort{\mathop{\hbox{\kern1pt\vrule width4.0pt height0.4pt depth0pt
                \vrule width0.4pt height6.0pt depth0pt\kern3.5pt}}}
\let\lra\longrightarrow
\def\vlra{\mathrel{\smash-}\joinrel\mathrel{\smash-}\joinrel%
\kern-2pt\longrightarrow}


\def\rk{\mathop{\rm rk}\nolimits}
\def\codim{\mathop{\rm codim}\nolimits}

\def\Sing{\mathop{\rm Sing}\nolimits}

\def\Hom{\mathop{\rm Hom}\nolimits}
\def\Ker{\mathop{\rm Ker}\nolimits}
\def\Im{\mathop{\rm Im}\nolimits}


\long\def\InsertFig#1 #2 #3 #4\EndFig{\par
\hbox{\hskip #1mm$\vbox to#2mm{\vfil\special{"
(/home/demailly/psinputs/grlib.ps) run
#3}}#4$}}
\long\def\LabelTeX#1 #2 #3\ELTX{\rlap{\kern#1mm\raise#2mm\hbox{#3}}}



\title{A Kawamata-Viehweg Vanishing Theorem}
\title{on compact K\"ahler manifolds}
\titlerunning{A Kawamata-Viehweg Vanishing Theorem on compact 
K\"ahler manifolds}
\vskip10pt

{\noindent\hangindent0.6cm\hangafter-1
{\bf Jean-Pierre Demailly${}^\star$, Thomas Peternell${}^{\star\star}$

}}

{\noindent\hangindent0.6cm\hangafter-4{\sl
\llap{${}^\star~$}Universit\'e de Grenoble I, BP 74\hfill
${}^{\star\star}~$Universit\"at Bayreuth\kern0.6cm\break
Institut Fourier, UMR 5582 du CNRS
\hfill Mathematisches Institut\kern0.6cm\break
38402 Saint-Martin d'H\`eres, France\hfill D-95440 Bayreuth,
Deutschland\kern0.6cm

}}
\vskip20pt

\noindent{\bf Abstract.} We prove a Kawamata-Viehweg vanishing theorem 
on a normal compact K\"ahler space $X$: if $L$ is a nef line bundle
with $L^2 \ne 0$, then $H^q(X,K_X+L) = 0$ for $q \geq \dim X - 1$. As
an application we complete a part of the abundance theorem for minimal
K\"ahler threefolds: if $X$ is a minimal K\"ahler threefold, then the
Kodaira dimension $\kappa(X)$ is nonnegative.
\vskip20pt

\section{\S0. Introduction}

In this paper we establish the following Kawamata-Viehweg type
vanishing theorem on a compact K\"ahler manifold or, more generally, a
normal compact K\"ahler space.

\claim 0.1\ Theorem| Let $X$ be a normal compact K\"ahler space of 
dimension $n$ and $L$ a nef line bundle on $X$. Assume that 
$L^2 \ne 0$. Then
$$
H^{q}(X,K_X+L) = 0
$$
for $q\geq n-1$.
\endclaim

In general, one expects a vanishing 
$$
H^q(X,K_X+L) = 0
$$
for $q \geq n+1-\nu(L)$, where $\nu(L)$ is the numerical Kodaira
dimension of the nef line bundle $L$, i.e.\ $\nu(L)$ is the largest
integer $\nu$ such that $L^{\nu} \ne 0$.\smallskip 

\noindent
Of course, when $X$ is projective, Theorem 0.1 is contained in the
usual Kawamata-Viehweg vanishing theorem, but the methods of proof in
the algebraic case clearly fail in the general K\"ahler setting.
Instead we proceed in the following way. Clearly we may assume that
$X$ is smooth and by Serre duality, only the cohomology group
$H^{n-1}$ is of interest. Take a singular metric $h$ on $L$ with
positive curvature current $T$ with local weight function $h$.  
By [Si74, De93a] there exists a decomposition
$$
T = \sum \lambda_j D_j + G, \eqno (D)
$$
where $\lambda_j \geq 1$ are irreducible divisors, and $G$ is a
pseudo-effective current such that $G \vert D_i$ is pseudo-effective
for all $i$.  Consider the multiplier ideal sheaf $\cI(h)$. We
associate to $h$ another, ``upper regularized'' multiplier ideal sheaf
$\cI_+(h)$ by setting
$$
\cI_+(h):=\lim_{\varepsilon\to 0_+}\cI(h^{1+\varepsilon})
=\lim_{\varepsilon\to 0_+}\cI\big((1+\varepsilon)\varphi\big).
$$
It is unknown whether $\cI(h)$ and $\cI_+(h)$
actually differ; in all known examples they are equal. Then in Section 2 
the following vanishing theorem is proved.

\claim 0.2\ Theorem| Let $(L,h)$ be a holomorphic line bundle over 
a compact K\"ahler \hbox{$n$-fold}~$X$. Assume that $L$ is nef and has
numerical Kodaira dimension $\nu(L)=\nu\ge 0$, i.e.\ $c_1(L)^\nu\ne 0$ and
$\nu$ is maximal. Then the morphism
$$
H^q(X,\cO(K_X+L)\otimes\cI_+(h))\la H^q(X,K_X+L)
$$
induced by the inclusion $\cI_+(h)\subset\cO_X$ vanishes for $q>n-\nu$.
\endclaim

The strategy of the proof of Theorem 0.2 is based on a direct
application of the Bochner technique with special hermitian metrics
constructed by means of the
Calabi-Yau theorem.\smallskip

\noindent
Now, coming back to the principles of the proof of Theorem 0.1, we
introduce the divisor
$$
D = \sum  [\lambda_j]D_j.
$$ 
Then Theorem 0.2 yields the vanishing of the map in cohomology
$$
H^{n-1}(X,-D+L+K_X) \la H^{n-1}(X,L+K_X).
$$
Thus we are reduced to show that $H^{n-1}(D,L+K_X \vert D) = 0$, or dually
that
$$
H^0(D,-L+D \vert D) = 0.
$$ 
This is now done by a detailed analysis of a potential non-zero section 
in $-L+D \vert D$; making use of the decomposition $(D)$ and of a Hodge 
index type inequality.\smallskip

\noindent
The vanishing theorem 0.1 is most powerful when $X$ is a threefold,
and in the second part of the paper we apply 0.1 - or rather a
technical generalization - to prove the following abundance theorem.

\claim 0.3\ Theorem| Let $X$ be a $\bQ$-Gorenstein K\"ahler threefold
with only terminal singularities, such that $K_X$ is nef (a minimal
K\"ahler threefold for short). Then $\kappa (X) \geq 0$.
\endclaim

This theorem was established in the projective case by Miyaoka and in
[Pe01] for K\"ahler threefolds, with the important exception
that $X$ is a simple threefold which is not Kummer. Recall that $X$ is
said to be {\it simple} if there is no proper compact subvariety
through a very general point of $X$, and that $X$ is said to be Kummer
if $X$ is bimeromorphic to a quotient of a torus. So our contribution
here consists in showing that such a simple threefold $X$ with $K_X$
nef has actually $\kappa (X) = 0$. Needless to say that among all
K\"ahler threefolds the simple non-Kummer ones (which conjecturally
do not exist) are most difficult to deal with, since they do not carry 
much global information besides the fact that $\pi_1$ is finite and 
that they have a holomorphic $2$-form.
\medskip
\noindent The first main ingredient in our approach is the inequality
$$
K_X \cdot c_2(X) \geq 0
$$
for a minimal simply connected K\"ahler threefold $X$ with algebraic
dimension $a(X) = 0$. Philosophically this inequality comes from
Enoki's theorem that the tangent sheaf of $X$ is $K_X$-semi-stable
when $K_X^2 \ne 0$ resp.\ $(K_X,\omega)$-semi-stable when $K_X^2 = 0;$
here $\omega$ is any K\"ahler form on $X$. Now if this semi-stability
with respect to a degenerate polarization would yield a Miyaoka-Yau
inequality, then $K_X \cdot c_2(X) \geq 0$ would follow. However this
type of Miyoka-Yau inequalities with respect to degenerate polarizations is
completey unknown. In the projective case, the inequality follows from
Miyaoka's generic nefness theorem and is based on char. $p$-methods.
Instead we approximate $K_X$ (in cohomology) by K\"ahler forms
$\omega_j$. If $T_X$ is still $\omega_j$-semi-stable for sufficiently
large $j$, then we can apply the usual Miyaoka-Yau inequality and pass
to the limit to obtain $K_X \cdot c_2(X) \geq 0$. Otherwise we examine
the maximal destabilizing subsheaf which essentially (because of $a(X)
= 0$) is independent of the polarization.

\smallskip \noindent
The second main ingredient is the boundedness $h^2(X,mK_X) \leq 1$. If
$K_X^2 \ne 0$, this is of course contained in Theorem 0.1. If $K_X^2 =
0$, we prove this boundedness under the additional assumption that
$a(X) = 0$ and that $\pi_1(X) $ is finite (otherwise by a result of
Campana $X$ is already Kummer). The main point is that if $h^2(X,mK_X)
\geq 2$, then we obtain ``many'' non-split extensions
$$ 0 \la K_X \la \cE \la \ mK_X \la 0$$
and we analyze whether $ \cE$ is semi-stable or not. The assumption on
$\pi_1$ is used to conclude that if $\cE$ is projectively flat, then
$\cE$ is trivial after a finite \'etale cover.\smallskip 

\noindent
From these two ingredients Theorem 0.3 immediately follows by
applying Riemann-Roch on a desingularization of $X$.

\bigskip \noindent
The only remaining problem concerning abundance on K\"ahler threefolds is
to prove that a simple K\"ahler threefold with $K_X$ nef and
$\kappa (X) = 0$ must be Kummer.
\vskip20pt

\section{\S1. Preliminaries}

We start with a few preliminary definitions.

\claim 1.1\ Definition| A normal complex space $X$ is said to be
K\"ahler if there exists a K\"ahler form $\omega$ on the regular part
of $X$ such that the following holds. Every singular point $x \in X$ admits
an open neighborhood $U$ and a closed embedding $U \subset V$ into an open
set $U \subset \bC ^N$ such that there is a K\"ahler form $\eta$ on
$V$ with $\eta \vert U = \omega$.
\endclaim

\claim 1.2\ Remark| {\rm  Let $X$ be a compact K\"ahler space and let
$f: \hat X \la X$ be a desingularization by a sequence of blow-ups.
Then $\hat X$ is a K\"ahler manifold. More generally consider a 
holomorphic map $f: \hat X \la X$ of a normal compact complex space
to a normal compact K\"ahler space. If $f$ is a projective morphism
or, more generally, a K\"ahler morphism, then $\hat X$ is K\"ahler.
For references to this and more informations on K\"ahler spaces, 
we refer to [Va84].
\smallskip
\noindent A K\"ahler form $\omega$ defines naturally a class $[\omega] \in 
H^2(X,\bR)$, see
[Gr62] where K\"ahler metrics on singular spaces were first introduced. 
Therefore we also have a K\"ahler cone on a normal variety.  }
\endclaim

\claim 1.3\ Notation| {\rm  Let $X$ be a normal compact complex space.
{{\item {\rm(1)} Let $A$ and $B$ be reflexive sheaves of rank $1$. Then
we define $A \hotimes B := (A \otimes B)^{**}$. Moreover we let
$A^{[m]} := A^{\hotimes m}.$}
{\item {\rm(2)} A reflexive sheaf $A$ is said to be a $\bQ$-line bundle 
if there exists a positive integer $m$ such that $A^{[m]}$ is locally
free. }
{\item {\rm(3)} $X$ is $\bQ$-Gorenstein if the canonical reflexive sheaf
$\omega_X$, also denoted $K_X$, is a $\bQ$-line bundle. $X$ is $\bQ$-factorial,
if every reflexive sheaf of rank $1$ is a $\bQ$-line bundle. }}}
\endclaim
\bigskip

\noindent

\claim 1.4\ Definition|  Let $X$ be a normal compact K\"ahler threefold.
{ {\item{\rm(1)} $X$ is simple if there is no proper compact subvariety 
through the very general point of $X.$}
{\item{\rm(2)} $X$ is Kummer, if $X$ is bimeromorphic to a quotient $T/G$
where $T$ is a torus and $G$ a finite group acting on $T.$}}
\endclaim

It is conjectured that all simple threefolds are Kummer. 

\claim 1.5\ Notation| {\rm {{\item{(1)} The {\it algebraic dimension} $a(X)$
of an irreducible reduced compact complex space is the transcendence
degree of the field of meromorphic functions over $\bC$. If $a(X) =
0$, i.e.\ all meromorphic functions on $X$ are constant, then it is
well known that $X$ carries only finitely many irreducible hypersurfaces.}  
{\item{(2)} A line bundle $L$ on a compact K\"ahler
manifold is {\it nef}, if $c_1(L) $ lies in the closure of the K\"ahler
cone. For alternative descriptions see e.g. [DPS94,00]. If $X$ is a
normal compact K\"ahler, then $L$ is nef if there exists a
desingularization $\pi: \hat X \la X$ such that $\pi^*(L)$ is nef. By
[Pa98], this definition does not depend on the choice of $\pi.$}}}
\endclaim \vskip20pt

\section{\S2. Hodge index type inequalities}

We give here some generalizations of Hodge index inequalities for
nef classes over compact K\"ahler manifolds. In this direction the main 
result is the Hovanskii-Teissier concavity inequality, which can 
be stated in the following way (see e.g.\ [De93b], Prop.\ 5.2 and 
Remark 5.3).

\claim 2.1\ Proposition| Let $\alpha_1,\ldots,\alpha_k$ and $\gamma_1,\ldots
\gamma_{n-k}$ be nef cohomology classes on a compact K\"ahler $n$-dimensional
manifold $X$. Then
$$
\alpha_1 \cdots \alpha_k\cdot \gamma_1\cdots \gamma_{n-k}\ge
(\alpha_1^k\cdot \gamma_1\cdots \gamma_{n-k})^{1/k}\cdots
(\alpha_k^k\cdot \gamma_1\cdots \gamma_{n-k})^{1/k}.
$$
\endclaim

We want to derive from these a non vanishing property for intersection
products of the form $\alpha^i\cdot\beta^j$. Let us fix 
a K\"ahler metric $\omega$ on $X$. By Proposition 2.1 applied with $k=i+j$ and
the $\alpha_\ell$'s being $i$ copies of $\alpha$ followed by $j$ copies of
$\beta$ and $\gamma_\ell=\omega$, we have
$$
\alpha^i\cdot\beta^j\cdot\omega^{n-i-j}\ge 
(\alpha^k\cdot \omega^{n-k})^{i/k}\cdots(\beta^k\cdot \omega^{n-k})^{j/k}.
$$
As all products $\alpha^k$ and analogues can be represented by closed positive
currents, we have $\alpha^k\ne 0\Rightarrow \alpha^k\cdot \omega^{n-k}>0$,
hence with $k=i+j$ we find
$$
\alpha^{i+j}\ne 0~~\hbox{and}~~\beta^{i+j}\ne 0~~\Longrightarrow~~
\alpha^i\cdot\beta^j\ne 0.\leqno(2.2)
$$
This is of course optimal in terms of the exponents if $\alpha=\beta$, but
as we shall see in a moment, this is possibly not optimal in a dissymetric 
situation.  Actually, we have the following additional inequalities
which can be viewed as ``differentiated'' Hovanskii-Teissier
inequalities.

\claim 2.3\ Theorem| Let $\alpha$ and $\beta$ be nef cohomology classes
of type $(1,1)$ on a compact K\"ahler $n$-dimensional manifold $X$. Assume 
that $\alpha^p\ne 0$ and $\beta^q\ne 0$ for some integers $p,q>0$. Then 
we have $\alpha^i\cdot \beta^j\ne 0$ as soon as there exists an
integer $k\ge i+j$ such that
$$i(k-p)_+ + j(k-q)_+<k,$$
where $x_+$ means the positive part of a number $x$.
\endclaim

\proof. Assume that $\alpha^i\cdot\beta^j=0$. We apply the 
Hovanskii-Teissier inequality respectively with 
$\alpha_\ell=\alpha+\varepsilon\omega$ ($i$ terms), or 
$\alpha_\ell=\beta+\varepsilon\omega$ ($j$ terms) or
$\alpha_\ell=\omega$ ($k-i-j$ terms), and $\gamma_\ell=\omega$. This gives
$$
(\alpha+\varepsilon\omega)^i\cdot(\beta+\varepsilon\omega)^j
\cdot\omega^{n-i-j}\ge
\big((\alpha+\varepsilon\omega)^k\cdot\omega^{n-k}\big)^{i/k}
\big((\beta+\varepsilon)^k\cdot\omega^{n-k}\big)^{j/k}
(\omega^n)^{1-i/k-j/k}.\leqno(*)
$$
By expanding the intersection form and using the assumption
$\alpha^i\cdot\beta^j=0$, we infer
$$
(\alpha+\varepsilon\omega)^i\cdot(\beta+\varepsilon\omega)^j
\cdot\omega^{n-i-j}\le O(\varepsilon)
$$
as $\varepsilon$ tends to zero. On the other hand
$(\alpha+\varepsilon\omega)^k\cdot\omega^{n-k}$ is bounded away from $0$ if
$k\le p$ since then $\alpha^k\ne 0$, and
$(\alpha+\varepsilon\omega)^k\cdot\omega^{n-k}\ge C\varepsilon^{k-p}$ for some
constant $C>0$ if $k\ge p$. Hence we infer from $(*)$ that
$$
C\varepsilon^{(i/k)(k-p)_++(j/k)(k-q)_+}=O(\varepsilon),
$$
and this is not possible if $i(k-p)_++j(k-q)_+<k$. The theorem is proved.\qed
\medskip

\noindent
The special case $p=2$, $q=1$, $i=j=1$, $k=2$ provides the following 
result which will be needed later on several occasions.

\claim 2.4\ Corollary| Assume that $\alpha$, $\beta$ are nef with 
$\alpha^2\ne 0$ and $\beta\ne 0$. Then $\alpha\cdot\beta\ne 0$.
\endclaim

Finally, we state an extension of Proposition 2.1 in the
case when one of the cohomology classes involved is not necessarily nef.

\claim 2.5\ Proposition|Let $\alpha$ be a real $(1,1)$-cohomology
class, and let $\beta$, $\gamma_1,\ldots \gamma_{n-2}$ be nef
cohomology classes. Then
$$
(\alpha \cdot \beta\cdot \gamma_1\cdots \gamma_{n-2})^2\ge
(\alpha^2\cdot \gamma_1\cdots \gamma_{n-2})
(\beta^2\cdot \gamma_1\cdots \gamma_{n-2}).
$$
\endclaim

\proof.. By proposition 2.1, the result is true when $\alpha$ is nef. If we
replace $\beta$ by $\beta+\varepsilon\omega$ and let $\varepsilon>0$ tend to
zero, we see that it is enough to consider the case when $\beta$ is a
K\"ahler class. Then $\alpha+\lambda\beta$ is also K\"ahler for $\lambda\gg 1$
large enough, and the inequality holds true with $\alpha+\lambda\beta$
in place of $\alpha$. However, after making the replacement, the contributions
of terms involving $\lambda$ in the right and left hand side of the inequality
are both equal to
$$
2\lambda(\alpha \cdot \beta\cdot \gamma_1\cdots \gamma_{n-2})
(\beta^2\cdot \gamma_1\cdots \gamma_{n-2})+
\lambda^2(\beta^2\cdot \gamma_1\cdots \gamma_{n-2})^2.
$$
Hence these terms cancel and the claim follows.\qed

\section{\S3. Partial vanishing for multiplier ideal sheaf cohomology}

Let $(L,h)$ be a holomorphic line bundle over a compact K\"ahler
$n$-fold $X$. Locally in a trivialization $\smash{L_{|U}}\simeq U\times\bC$,
the metric is given by $\Vert\xi\Vert_x=|\xi|e^{-\varphi(x)}$ and
we assume that the curvature $\Theta_h(L):={i\over\pi}\ddbar\varphi$
is a closed positive current (so that, in particular, $L$ is
pseudo-effective). We introduce as usual the {\it multiplier ideal sheaf}
$\cI(h):=\cI(\varphi)$ where
$$
\cI(\varphi)_x:=\big\{
f\in\cO_{X,x}\,;\;\exists V\ni x,~\int_V|f(z)|^2e^{-2\varphi(z)}<+\infty\big\}
$$
and $V$ is an arbitrarily small neighborhood of~$x$. We also consider 
the {\it upper regularized multiplier ideal sheaf}
$$
\cI_+(h):=\lim_{\varepsilon\to 0_+}\cI(h^{1+\varepsilon})
=\lim_{\varepsilon\to 0_+}\cI\big((1+\varepsilon)\varphi\big).
$$
It should be noticed that $\cI\big((1+\varepsilon)\varphi\big)$ increases as
$\varepsilon$ decreases, hence the limit is locally stationary by the
Noether property of coherent sheaves, and one has of course
\hbox{$\cI_+(h)\subset\cI(h)$}. It is unknown whether these sheaves may
actually differ (in all known examples they are equal). In any case, they
coincide at least in codimension $1$ (i.e., outside an analytic subset 
of codimension${}\ge 2$).

\claim 3.2\ Proposition|Let
$$
\Theta_h(L)=\sum_{j=1}^{+\infty}\lambda_jD_j+G
$$
be the Siu decomposition of the $(1,1)$-current $\Theta_h(L)$ as a countable
sum of effective divisors and of a $(1,1)$-current $G$ such that the Lelong
sublevel sets $E_c(G)$, $c>0$, all have codimension $2$. Then we have the
inclusion of sheaves
$$
\cI_+(h)\subset\cI(h)\subset\cO\big(-\sum[\lambda_j]D_j\big),\qquad
[\lambda_j]:=\hbox{integer part of $\lambda_j$},
$$
and equality holds on $X\ssm Z$ where $Z$ is an analytic subset of $X$ whose
components all have codimension at least $2$.
\endclaim

\proof.. The decomposition exists by [Siu74] (see also [De93a]). Now, if $g_j$
is a local generator of the ideal sheaf $\cO(-D_j)$, the plurisubharmonic
weight $\varphi$ of $h$ can be written as
$$
\varphi=\sum\lambda_j\log|g_j|+\psi
$$
where $\psi$ is plurisubharmonic and the $E_c(\psi)$ have codimension $2$ 
at least. Since $\psi$ is locally bounded from above, it is obvious
that 
$$\cI(\varphi)\subset\cI(\lambda_j\log|g_j|)\subset
\cO\big(-\sum[\lambda_j]D_j\big).$$
Now, let $Y$ be the union of all sets $E_c(\psi)$ 
(with, say, $c=1/k$), all pairwise
intersections $D_j\cap D_k$ and all singular sets $D_{j\,{\rm sing}}$. This
set $Y$ is at most a countable union of analytic sets of codimension${}\ge 2$.
Pick an arbitrary point $x\in X\ssm Y$. Then $x$ meets the support of
$\bigcup D_j$ in at most one point which is then a smooth point of some 
$D_k$, and the Lelong number of $\psi_k=\psi+\sum_{j\ne k}
\lambda_j\log|g_j|$ at $x$ is zero. Then $\varphi=\lambda_k\log|g_k|+\psi_k$
and the inclusion
$$
\cI(h)_x \supset \cO\big(-\sum_j[\lambda_j]D_j\big)_x=
\cO(-[\lambda_k]D_k)_x
$$
holds true by H\"older's inequality. In fact, for every germ $f$ in
$\cO(-[\lambda_k]D_k)_x$ we have
$$
\int_{V\ni x}|f|^2\exp\big(-(1+\varepsilon)\lambda_k\log|g_k|\big)<+\infty
$$
for $\varepsilon>0$ so small that $[(1+\varepsilon)\lambda_k]=[\lambda_k]$,
while $e^{-\psi_k}$ is in $L^p(V_p)$ for some $V_p\ni x$, for every $p>1$. 
Similarly, we have
$$
\cI_+(h)_x \supset\cO\big(-[(1+\varepsilon)\lambda_k]D_k\big)_x=
\cO\big(-[\lambda_k]D_k\big)_x
$$
for $\varepsilon>0$ small enough. The analytic set $Z$ where our sheaves 
differ $\big[$i.e.\ the union of supports of $\cI(h)/\cI_+(h)$ and
$\cO\big(-\sum_j[\lambda_j]D_j\big)/\cI(h)\big]$ must be contained in 
$Y$, hence $Z$ is of codimension${}\ge 2$.\qed
\medskip

The main goal of this section is to prove the following partial
vanishing theorem.

\claim 3.3\ Theorem| Let $(L,h)$ be a holomorphic line bundle over 
a compact K\"ahler \hbox{$n$-fold}~$X$, equipped with a singular metric $h$
such that $\Theta_h(L)\ge 0$. Assume that $L$ is nef and has
numerical dimension $\nu(L)=\nu\ge 0$, i.e.\ $c_1(L)^\nu\ne 0$ and
$\nu$ is maximal. Then the morphism
$$
H^q(X,\cO(K_X+L)\otimes\cI_+(h))\la H^q(X,K_X+L)
$$
induced by the inclusion $\cI_+(h)\subset\cO_X$ vanishes for $q>n-\nu$.
\endclaim

Of course, it is expected that the Kawamata-Viehweg vanishing theorem 
also holds for K\"ahler manifolds, in which case the whole group 
$H^q(X,K_X+L)$ vanishes and Theorem 3.3 would then be an obvious consequence. 
However, we will see in Section~4 that, conversely, Theorem 3.3 can be 
used to derive the Kawamata-Viehweg vanishing theorem in the first 
non trivial case $\nu=2$. Using the same method for higher values of $\nu$
would probably be very hard, if not impossible.

\proof.. Our strategy is based on a direct application of the Bochner 
technique with special hermitian metrics constructed by means of the 
Calabi-Yau theorem.

Let us fix a smooth hermitian metric $h_\infty$ on $L$, which may have
a curvature form $\Theta_{h_\infty}(L)$ of arbitrary sign, and let 
$\varepsilon>0$. Then $c_1(L)+\varepsilon\omega$ is a K\"ahler
class, hence by the Calabi-Yau theorem for complex Monge-Amp\`ere equations
there exists a hermitian metric 
$h_\varepsilon=h_\infty e^{-2\varphi_\varepsilon}$ such that
$$
\big(\Theta_{h_\varepsilon}(L)+\varepsilon\omega\big)^n=C_\varepsilon\omega^n.
\leqno(3.4)
$$
Here $C_\varepsilon>0$ is the constant such that
$$
C_\varepsilon={\int_X(c_1(L)+\varepsilon\omega)^n\over\int_X\omega^n}\ge
C\varepsilon^{n-\nu}.
$$
Let $h=h_\infty e^{-2\psi}$ be a metric with $\Theta_h(L)\ge 0$ as given
in the statement of the theorem, and let $\psi_\varepsilon\downarrow\psi$ 
be a regularization of $\psi$ possessing only analytic singularities 
(i.e.\ only logarithmic poles), such that 
$$
\tilde h_\varepsilon:=h_\infty e^{-2\psi_\varepsilon}
$$
satisfies $\Theta_{\tilde h_\varepsilon}(L)\ge-\varepsilon\omega$ in the
sense of currents. Such a metric exists by the general regularization results
proved in [De92]. We consider the metric
$$
\hat h_\varepsilon=(h_\varepsilon)^\delta(\tilde h_\varepsilon)^{1-\delta}=
h_\infty\exp\big(-2(\delta\varphi_\varepsilon+(1-\delta)\psi_\varepsilon)\big)
$$
where $\delta>0$ is a sufficiently small number which will be fixed later. By
construction,
$$
\eqalign{
\Theta_{\hat h_\varepsilon}(L)+2\varepsilon\omega&=
\delta\big(\Theta_{h_\varepsilon}(L)+\varepsilon\omega\big)+
(1-\delta)\big(\Theta_{\tilde h_\varepsilon}(L)+\varepsilon\omega\big)+
\varepsilon\omega\cr
&\ge
\delta\big(\Theta_{h_\varepsilon}(L)+\varepsilon\omega\big)+\varepsilon\omega.
\cr}
$$
Denote by $\lambda_1\le\ldots\le\lambda_n$ and
$\hat\lambda_1\le\ldots\le\hat\lambda_n$, respectively, the eigenvalues of 
the curvature forms $\Theta_{h_\varepsilon}(L)+\varepsilon\omega$ and
$\Theta_{\hat h_\varepsilon}(L)+2\varepsilon\omega$ at every point $z\in X$, 
with respect to the base K\"ahler metric $\omega(z)$. By the minimax 
principle we find $\hat\lambda_j\ge\delta\lambda_j+\varepsilon$. On the other
hand, the Monge-Amp\`ere equation (3.4) tells us that
$$
\lambda_1\ldots\lambda_n=C_\varepsilon\ge C\varepsilon^{n-\nu}
\leqno(3.5)
$$
everywhere on $X$. We apply the basic Bochner-Kodaira inequality to
sections of type $(n,q)$ with values in the hermitian line bundle
$(L,\hat h_\varepsilon)$. As the curvature eigenvalues of
$\Theta_{\hat h_\varepsilon}(L)$ are equal to 
$\hat \lambda_j-2\varepsilon$ by definition, we find
$$
\Vert\dbar u\Vert^2_{\hat h_\varepsilon}+
\Vert\dbar^\star u\Vert^2_{\hat h_\varepsilon}\ge
\int_X(\hat\lambda_1+\cdots+\hat\lambda_q-2q\varepsilon)
|u|^2_{\hat h_\varepsilon}dV_\omega\leqno(3.6)
$$ 
for every smooth $(n,q)$-form $u$ with values in $L$. Actually this is
formally true only if the metric $\hat h_\varepsilon$ is smooth on $X$. 
The metric $h_\varepsilon$ is indeed smooth, but $\tilde h_\varepsilon$ 
may have poles along an analytic set $Z_\varepsilon\subset X$. In that case,
we apply instead the inequality to forms $u$ which are compactly
supported in $X\ssm Z_\varepsilon$, and replace the K\"ahler metric
$\omega$ by a sequence of complete K\"ahler metrics $\omega_k\downarrow
\omega$ on $X\ssm Z_\varepsilon$, and pass to the limit as $k$ tends to
$+\infty$ (see e.g.\ [De82] for details about such techniques). In the 
limit we recover the same estimates as if we were in the smooth case, 
and we therefore allow ourselves to ignore this minor technical
problem from now on.

Now, let us take a cohomology class $\{\beta\}\in H^q(X,K_X\otimes L\otimes
\cI_+(h))$. By using \v Cech cohomology and the De Rham-Weil isomorphism
between \v Cech and Dolbeault cohomology (via a partition of unity and the
usual homotopy formulas), we obtain a representative $\beta$ of the 
cohomology class which is a smooth $(n,q)$-form with values
in $L$, such that the coefficients of $\beta$ lie in the sheaf
$\cI_+(h)\otimes_{\cO_X}\cC^\infty$. We want to show that $\beta$ is a
boundary with respect to the cohomology group $H^q(X,K_X\otimes L)$.
This group is a finite dimensional Hausdorff vector space whose topology is 
induced by the $L^2$ Hilbert space topology on the space of forms 
(all Sobolev norms induce in fact the same topology on the level of 
cohomology groups). Therefore, it is enough to
show that we can approach $\beta$ by $\dbar$-exact forms in $L^2$ norm.

As in H\"ormander [H\"o65], we write every form $u$ in the domain of 
the $L^2$-extension of $\dbar^*$ as $u=u_1+u_2$ with 
$$
u_1\in\Ker\dbar\quad\hbox{and}\quad u_2\in(\Ker\dbar)^\perp=\overline{
\Im\dbar^*}\subset\Ker\dbar^*.
$$
Therefore, since $\beta\in\Ker\dbar$,
$$
\big|\Langle \beta, u\Rangle\big|^2=
\big|\Langle \beta, u_1\Rangle\big|^2\le
\int_X{1\over\hat\lambda_1+\cdots+\hat\lambda_q}|\beta|^2_{\hat h_\varepsilon}
dV_\omega
\int_X(\hat\lambda_1+\cdots+\hat\lambda_q)|u_1|^2_{\hat h_\varepsilon}
dV_\omega.
$$
As $\dbar u_1=0$, an application of (3.6) to $u_1$ (together with an 
approximation of $u_1$ by compactly supported smooth sections on the
corresponding complete K\"ahler manifold $X\ssm Z_\varepsilon$) shows that the 
second integral in the right hand side is bounded above by
$$
\Vert\dbar^* u_1\Vert_{\hat h_\varepsilon}^2+2q\varepsilon
\Vert u_1\Vert_{\hat h_\varepsilon}^2\le
\Vert\dbar^* u\Vert_{\hat h_\varepsilon}^2+2q\varepsilon
\Vert u\Vert_{\hat h_\varepsilon}^2,
$$
so we finally get
$$
\big|\Langle \beta, u\Rangle\big|^2\le
\int_X{1\over\hat\lambda_1+\cdots+\hat\lambda_q}|\beta|^2_{\hat h_\varepsilon}
dV_\omega
\big(\Vert\dbar^* u\Vert_{\hat h_\varepsilon}^2+2q\varepsilon
\Vert u\Vert_{\hat h_\varepsilon}^2\big).
$$
By the Hahn-Banach theorem (or rather a Hilbert duality argument in this
situation), we can find elements $v_\varepsilon$, $w_\varepsilon$ such that
$$
\Langle \beta, u\Rangle=
\Langle v_\varepsilon, \dbar^* u\Rangle+
\Langle w_\varepsilon, u\Rangle\quad\forall u,\qquad\hbox{i.e.}\quad
\beta=\dbar v_\varepsilon + w_\varepsilon,
$$
with
$$
\Vert v_\varepsilon\Vert_{\hat h_\varepsilon}^2+
{1\over 2q\varepsilon}\Vert w_\varepsilon\Vert_{\hat h_\varepsilon}^2\le
\int_X{1\over\hat\lambda_1+\cdots+\hat\lambda_q}|\beta|^2_{\hat h_\varepsilon}
dV_\omega.
$$
As a consequence, the $L^2$ distance of $\beta$ to the space of $\dbar$-exact
forms is bounded by $\Vert w_\varepsilon\Vert_{\hat h_\varepsilon}$ where
$$
\Vert w_\varepsilon\Vert_{\hat h_\varepsilon}^2=
\int_X|w_\varepsilon|^2_{h_\infty}e^{-2(\delta\varphi_\varepsilon+
(1-\delta)\psi_\varepsilon)}dV_\omega\le 2q\varepsilon
\int_X{1\over\hat\lambda_1+\cdots+\hat\lambda_q}|\beta|^2_{\hat h_\varepsilon}
dV_\omega.
$$
We normalize the choice of the potentials $\varphi_\varepsilon$, $\psi$ and
$\psi_\varepsilon$ so that
$$\sup_X\varphi_\varepsilon=0,\qquad
\sup_X\psi=-1,\qquad -1\le \sup_X\psi_\varepsilon<0\;;$$
in this way $\varphi_\varepsilon,\,\psi_\varepsilon\le 0$
everywhere on $X$ (all inequalities can be achieved simply by adding 
suitable constants). From this we infer
$$
\int_X|w_\varepsilon|^2_{h_\infty}dV_\omega\le 2\int_X{q\varepsilon
\over\hat\lambda_1+\cdots+\hat\lambda_q}|\beta|^2_{\hat h_\varepsilon}
dV_\omega,
$$
and what remains to be shown is that the right hand side converges to $0$
for a suitable choice of $\delta>0$. By construction 
$\hat\lambda_j\ge\delta\lambda_j+\varepsilon$ and (3.5) implies
$$
\lambda_q^q\lambda_{q+1}\ldots\lambda_n\ge
\lambda_1\ldots\lambda_n\ge C\varepsilon^{n-\nu},
$$
hence
$$
{1\over\lambda_1+\cdots+\lambda_q}\le {1\over\lambda_q}\le C^{-1/q}
\varepsilon^{-(n-\nu)/q}(\lambda_{q+1}\ldots\lambda_n)^{1/q}.
$$
We infer
$$
\gamma_\varepsilon:={q\varepsilon\over\hat\lambda_1+\cdots+\hat\lambda_q}\le
\min\Big(1,{q\varepsilon\over\delta\lambda_q}\Big)\le
\min\big(1,C\delta^{-1}\varepsilon^{1-(n-\nu)/q}
(\lambda_{q+1}\ldots\lambda_n)^{1/q}\big).
$$
We notice that
$$
\int_X \lambda_{q+1}\ldots\lambda_n\, dV_\omega\le \int_X
(\Theta_{h_\varepsilon}(L)+\varepsilon\omega)^{n-q}\wedge\omega^q=
(c_1(L)+\varepsilon\{\omega\})^{n-q}\{\omega\}^q\le C'',
$$
hence the functions $(\lambda_{q+1}\ldots\lambda_n)^{1/q}$ are uniformly
bounded in $L^1$ norm as $\varepsilon$ tends to zero. Since
$1-(n-\nu)/q>0$ by hypothesis, we conclude that $\gamma_\varepsilon$
converges almost everywhere to $0$ as $\varepsilon$ tends to zero.
On the other hand
$$
|\beta|^2_{\hat h_\varepsilon}=|\beta|^2_{h_\infty}
e^{-2(\delta\varphi_\varepsilon+(1-\delta)\psi_\varepsilon)}
\le|\beta|^2_{h_\infty}e^{-2\delta\varphi_\varepsilon}e^{-2\psi}.
$$
Our assumption that the coefficients of $\beta$ lie in $\cI_+(h)$
implies that there exists $p'>1$ such that 
$\int_X|\beta|^2_{h_\infty}e^{-2p'\psi}dV_\omega<\infty$. Let
$p\in{}]1,+\infty[$ be the conjugate exponent such that ${1\over p}+
{1\over p'}=1$. By H\"older's inequality, we have
$$
\int_X\gamma_\varepsilon|\beta|^2_{\hat h_\varepsilon}dV_\omega\le
\Big(\int_X|\beta|^2_{h_\infty}e^{-2p\delta\varphi_\varepsilon}
dV_\omega\Big)^{1/p}
\Big(\int_X\gamma_\varepsilon^{p'}|\beta|^2_{h_\infty}
e^{-2p'\psi}dV_\omega\Big)^{1/p'}.
$$
As $\gamma_\varepsilon\le 1$, the Lebesgue dominated convergence theorem
shows that
$$\int_X\gamma_\varepsilon^{p'}|\beta|^2_{h_\infty}e^{-2p'\psi}dV_\omega$$
converges to $0$ as $\varepsilon$ tends to $0$. However, the family of
quasi plurisubharmonic functions $(\varphi_\varepsilon)$ is a bounded 
family with respect to the $L^1$ topology on the space of 
(quasi)-plurisubharmonic functions -- we use here the fact that the currents
$$\Theta_{\tilde h_\varepsilon}(L)=\Theta_{h_\infty}(L)+{i\over\pi}\ddbar
\varphi_\varepsilon\ge 0$$ 
all sit in the same cohomology
class; the boundedness of their normalized potentials then results from 
the continuity properties of the Green operator.
By standard results of complex potential theory, we conclude that 
there exists a small constant $\eta>0$ such that
$\int_Xe^{-2\eta\varphi_\varepsilon}dV_\omega$ is uniformly bounded.
By choosing $\delta\le\eta/p$, the integral
$\int_X|\beta|^2_{h_\infty}e^{-2p\delta\varphi_\varepsilon}dV_\varepsilon$
remains bounded and we are done.\qed

\section{\S4. Kawamata-Viehweg vanishing theorem for line
bundles of numerical dimension 2}

In this section we prove the Kawamata-Viehweg vanishing theorem for
the cohomology group of degree $n-1$ of nef line bundles $L$ with $L^2 \ne
0$ on compact K\"ahler spaces of dimension $n$. Furthermore we will
prove an extended version where $L$ can be a reflexive sheaf. This 
will be needed for proving the abundance theorem for K\"ahler threefolds.

\claim 4.1\ Theorem| Let $X$ be a normal compact K\"ahler space of 
dimension $n$ and $L$ a nef line bundle on $X$. Assume that 
$L^2 \ne 0$. Then
$$ H^{q}(X,K_X+L) = 0$$
for $q  \geq n-1$.
\endclaim

\proof. 
In a first step we reduce the proof to the case of a smooth space $X$ (this 
is comfortable but not really necessary; all arguments would also work 
in the singular setting as well). In fact, let
$\pi: \hat X \la X$ be a K\"ahler desingularization. Then, assuming our
claim in the smooth case, we have
$$ H^{q}(\hat X,K_{\hat X}+\pi^*(L)) = 0.$$
By the projection formula and the Grauert-Riemenschneider vanishing theorem
$$R^j\pi_*(K_{\hat X}) = 0,$$ it follows
$$ H^{q}(X,\pi_*(K_{\hat X}) \otimes L) = 0.$$
Since $\pi_*(K_{\hat X}) \subset K_X$ with cokernel supported in 
codimension at least $2$, namely on the singular locus of $X$, the
vanishing claim follows. 

\bigskip

\noindent So from now on, we assume $X$ smooth. In the case $q = n$, we have
$H^n(X,K_X+L)=H^0(X,-L)^*$ by Serre duality, and for $L$ nef, $-L$ has
no section unless $L$ is trivial. Therefore the only interesting case
is $q = n-1$. We introduce a singular metric
$h$ on $L$ with positive curvature current $T$. 
By [Siu74] and [De92, De93a] we obtain a decomposition
$$ T = \sum \lambda_j D_j + G,$$
where $\lambda_j \geq 1$ are irreducible divisors, and $G$ is a 
positive current such that $G$ has Lelong numbers in codimension${}\ge 2$
only -- so that in particular $G\vert D_i$ is pseudo-effective for all $i$.
Consider the multiplier ideal sheaf $\cI(h)$. By Proposition 3.2 we have
$$\cI_+(h) \subset \cI(h) \subset \cO_X(- \sum [\lambda_j]D_j)$$
with equality in codimension 1. We put
$$ D = \sum  [\lambda_j]D_j.$$

We consider the canonical map in cohomology$$ H^{n-1}(X,-D+L+K_X) \lra 
H^{n-1}(X,L+K_X) $$
which is vanishing by (3.3). In order to prove our claim it is
therefore sufficient to prove
$$ H^{n-1}(D,L+K_X \vert D) = 0.$$
By Serre duality and the adjunction formula, this comes down to show
$$ H^0(D,-L+D \vert D) = 0.$$ 
Supposing the contrary, we fix a non-zero section
$$ \sigma \in H^0(D,-L+D).$$
We choose $p_1,\ldots, p_k$ maximal so that 
$$ \sigma \in H^0(D,-L+D-\sum p_jD_j),$$
i.e.\ we choose $\tilde D = \sum p_j D_j \subset D$ maximal
such that $\sigma \vert \tilde D = 0$. $\big[$In this notation,
we view $\tilde D$ as the subscheme of $X$ defined by the
structure sheaf $\cO_X/\cO_X(-\tilde D)\big]$.
\smallskip
\noindent
Then $0 \leq p_i \leq [\lambda_i] $ for all $i \in I$, not all
$p_i = [\lambda_i]$, and we shall always consider $\sigma$ as a section of
$-L+D-\sum p_j D_j$. Denote
$$ c_i = {{\{\lambda_i\} + p_i }\over {\lambda_i}}.$$
Then we have $0 \leq c_i \leq 1$. We introduce $c=\min c_i$ and
$$ I_0 = \{ i \in I\,\vert\; c_i = c\}.$$
Clearly $c < 1$. Notice that by construction $\sigma \vert D_i \ne 0 $
unless $c_i = 1$. Let
$$ E = - \big (\sum ( \{\lambda_i\}+ p_i)D_i \big ) - G.$$
Since $L = \sum \lambda_i D_i + G$, we have
$$ -L+D - \sum p_i D_i = - \big ( \sum (\{\lambda_i\} + p_i)D_i \big ) - G 
= E,$$
so $E$ is effective (possibly zero) on every $D_i$ with $c_i < 1$.
Since $L$ is nef, also the $\bR$-divisor $cL = \sum \lambda_i c D_i + cG$
is nef. Adding this to the divisor $E$ in the last equation, we deduce that
$$ - \big (\sum ( \{\lambda_i \} + p_i - c \lambda_i) D_i \big ) - (1-c)G $$
is pseudo-effective on every $D_i$ with $c_i < 1$.
Since $\{\lambda_i\} + p_i - c \lambda_i = 0$ for $ i \in I_0$, it follows
that
$$ - \big (\sum_{i \not \in I_0}(\{\lambda_i\} + p_i - c \lambda_i)D_i \big ) -
(1-c)G$$
is pseudo-effective on every $ D_i$ with $c_i < 1$, in particular for
every $i\in I_0$.
Now $ D_j \vert D_i$ is effective (possibly $0$) for all $j \ne i$,
and $G \vert  D_i$ is always pseudo-effective, hence, having in mind 
$c < 1$ and $\{\lambda_i \} + p_i - c \lambda_i > 0$ for $i \not \in I_0$,
we conclude that
$$ D_j \vert D_i \equiv 0  \eqno (1)$$
for all $(j,i)$ with $j \not \in I_0$ and $i \in I_0$     and that
$$ G \vert D_i \equiv 0 \eqno (2) $$
for all $i \in I_0$. 
Introducing  $$D' = \sum_{i \in I_0} \lambda_i D_i $$ and
$$ D'' = \sum_{i \not \in I_0}\lambda_i D_i,$$
we have 
$$L = D' + D'' + G$$ 
and $D'' \cdot D_i = G \cdot D_i = 0$
for all $i \in I_0$ by (1) and (2).
Hence $L \cdot D_i = D' \cdot D_i$ for $i \in I_0$, so that $D' \vert D_i$
is nef, hence $D'$ is nef by [Pa98]. In total
$$ L \cdot D' = D' \cdot D' $$
and
$$ \ D' \cdot D'' = D' \cdot G = 0.$$
As $L^2\ne 0$ et $D' \ne 0$, Corollary 2.4 implies $L\cdot D'\ne 0$.
First recall that 
$$
L\vert D_i = \sum_{j \in I_0} \lambda_j D_j \vert D_i
$$
is nef. On the other hand 
$$-\sum_{j \in I_0} \{\lambda_j \} + p_j)D_j = 
- c \sum_{j \in I_0}\lambda_j D_j$$
is of course pseudo-effective on every $D_i$ for $i \in I_0$ ($E$ is effective
on those $D_i)$.
Combining these two facts, we deduce that either $c = 0$ or that 
$L \cdot D_i = 0$ for all $i \in I_0$, hence $L \cdot D' = 0$,
contradiction. So we have $c = 0$. This means $p_j = 0$ and $\lambda_j \in
\bN$ for all $j \in I_0$. 

\claim 4.2 Claim| The divisor $D''+G$ is nef, and in fact must be 
equal to zero.
\endclaim

\proof.\ {\it of the claim}. We consider the closed positive $(1,1)$-current
$\Theta=[D'']+G$. By the results of [Pa98], the proof of nefness of
$\{\Theta\}$ just amounts to showing that the restriction 
$\smash{\{\Theta\}_{|Z}}$ of the 
$(1,1)$-cohomology class $\{\Theta\}$ to any component $Z$ in the 
Lelong sublevel set $\bigcup_{c>0}E_c(\Theta)$ is nef. However 
$Z$ is either a component of $D''$ or a component of $\bigcup_{c>0}E_c(G)$. 
In the first case, $Z$ is contained in the support of $D''$, and as 
$D'\cdot D''=0$, $Z$ must be disjoint from $D'$. Hence
$$
\{\Theta\}_{|Z}=\{D''+G\}_{|Z}=\{L-D'\}_{|Z}=\{L\}_{|Z}
$$
is nef. If $Z$ is a component of $\bigcup_{c>0}E_c(G)$, then $Z$ has 
codimension at least $2$. Then we know by [De93a] that the intersection
product $[D']\wedge G$ is well defined as a closed positive current.
Since the cohomology class of this current is zero, we must
have $[D']\wedge G=0$. However, we infer from [De93a] that
$$\nu([D']\wedge G, z)\ge \nu([D'],z)\,\nu(G,z)>0$$
at every point $z\in D'\cap Z$, hence $Z$ must also be disjoint from 
$D'$ in that case. We conclude as
before that $\smash{\{\Theta\}_{|Z}}=\smash{\{L\}_{|Z}}$ is nef.
Now we have $D' \cdot (D''+G) = 0$, with $D'$, $D''+G$ nef and 
$D^{\prime 2}=L\cdot D'\ne 0$. Hence $\{D''+G\}=0$ by Corollary 2.4,
and we conclude that $D''=0$, $G=0$ (both $[D'']$ and $[G]$ being
positive currents).\qed
\medskip \noindent
From this we infer $L\equiv D'$ and $\cI(h) = \cO_X(-D')$. 

\medskip \noindent
{\bf Case 1} We assume that $L = D'$. 
Now the sequence 
$$
0 \la \cI(h) \otimes K_X + L \la K_X + L \la K_X + L \vert D' = K_{D'} \la 0
$$
gives in cohomology
$$
0 \la H^{n-1}(X,K_X + L) \la H^{n-1}(D',K_{D'}) \simeq 
H^0(D',\cO_{D'}) \la H^n(K_X) =
\bC \la 0.
$$ 
Thus we need to show
$$
h^0(D', \cO_{D'}) = 1.
$$ 
In order to verify this, we first observe that $D'$ is connected. In fact
otherwise write $D' = A+B$ with $A$ and $B$ effective and $A \cdot B = 0$.
But $A$ and $B$ are necessarily nef, hence the Hodge Index Theorem gives
a contradiction to $L^2 = (D')^2 \ne 0$. So $D'$ is connected and if 
$h^0(D',\cO_{D'}) \geq 2$, then $\cO_{D'}$ contains a nilpotent section $t
\ne 0$. 
Let $\sum_{j \in I}\mu_jD_j$ denote its vanishing divisor (notice that
$D'$ is Cohen-Macaulay!). Then $1 \leq \mu_j \leq \lambda_j$ for all
$j$. Let 
$$ J = \{ j \in J \vert {{\lambda_j} \over {\mu_j}} {\rm maximal}\}$$
and let $c = {{\lambda_j} \over {\mu_j}}$ be the maximal value. 
Notice that $- \sum_{j \in I}\mu_j D_j \vert D_i$ is effective
(possibly $0$) for all $i$. 
First we rule out the case that $c =   {{\lambda_j} \over {\mu_j}}$
for all $j \in I$. In fact, then $ L \vert D_i =c \sum \mu_j D_j \vert D_i $ 
is nef and its dual is effective, hence $L \vert D_i \equiv 0$ for all $i$,
whence $L^2 = 0$, contradiction.
Thus we find some $j$ such that 
$$
c>{{\lambda_j} \over {\mu_j}}.
$$ 
By connectedness of $D = D'$ we can choose $i_0 \in J$ in such a way that 
there exists $j_1 \in I \setminus J$ with $D_{i_0} \cap D_{j_1} \ne
\emptyset$. 
Now 
$$
\sum_{j \in I}(\lambda_j - c \mu_j)D_j \vert D_{i_0}
$$
is pseudo-effective as a sum of a nef and an effective line bundle
(this has nothing to do with the choice of $i_0$). Since the sum, taken over
$I$, is the same as the sum taken over $I \setminus \{i_0\}$, we conclude
that 
$$
\sum_{j \ne i_0}(\lambda_j - \mu_j)D_j \vert D_{i_0}
$$
is pseudo-effective, too. 
Now all $\lambda_j- c \mu_j \leq 0$ and $\lambda_{
j_1} - c \mu_{j_1} < 0$ with $D_{j_1} \cap D_{i_0} \ne \emptyset$,
hence the dual of
$$
\sum_{j \ne i_0}(\lambda_j - \mu_j)D_j \vert D_{i_0}
$$
is effective non-zero, a contradiction.\qed
\medskip 

\noindent
{\bf Case 2.} Now we  deal with the case that $L \not = D'$. Then we can 
write 
$$
L = D' + L_0
$$
where $L_0^m \in {\rm Pic}^0(X)$ (The exponent $m$ is there because there
might be torsion in $H^2(X,\bZ)$; we take $m$ to kill the denominator
of the torsion part). We may in fact assume that $m = 1$; otherwise we 
pass to a finite \'etale cover $\tilde X$ of $X$ and argue there
(the vanishing on $\tilde X$ clearly implies the vanishing on $X$). 
Then the sequence $S$ is modified to
$$ 0 \la \cI(h) \otimes (K_X + L) \la K_X + L \la (K_X + L)\vert D' 
= (K_{D'} + L_0) \vert D' \la 0. \leqno (S) $$ 
Taking cohomology as before, things come down to prove 
$$
H^0(D',-L_0 \vert D') = 0. \leqno (*)
$$
If $-L_0 \vert D' \ne 0$, then we see as above that $-L_0 \vert D'$
cannot have a nilpotent section. So if $(*)$ fails, then $-L_0 \vert D'$
has a section $s$ such that $s \vert {\rm red}D'$ has no zeroes, so
that $-L_0 \vert {\rm red}D'$ is trivial. But then $-L_0 \vert D'$ is
trivial. Now let $\alpha: X \la A$ be the Albanese map with image $Y$.
Then $L_0 = \alpha^*(L_0')$ with some line bundle $L_0'$ on $A$ which
is topologically trivial but not trivial.  Since $L_0 \vert D'$ is
trivial, we conclude that $\alpha (D') \ne Y$, and $\alpha (D') $ is
contained in a proper subtorus $B$ of $A$. Now consider the induced
map
$$
\beta: X \la A/B
$$
and denote its image by $Z$. Then $\beta (D')$ is a point; on the
other hand $D'$ is nef, so that $\dim Z = 1$ and $D'$ consists of
multiples of fibers of $\beta$. But this contradicts $D'^2 = L^2\ne 0$.\qed

For applications to minimal K\"ahler 3-folds, 4.1 is still 
not good enough, because we need to know the vanishing property 
$H^2(X,mK_X) = 0$ on a $\bQ$-Gorenstein $3$-fold (with $K_X^2 \ne 0$). 
We would like to set $L = (m-1)K_X$ to apply 4.1 but this is no
longer a line bundle. This difficulty is overcome by

\claim 4.3 Proposition| Let $X$ be a normal $\bQ$-Gorenstein compact 
K\"ahler 3-fold with at most terminal singularities. Let $A$ be a 
$\bQ$-line bundle. Suppose $A$ is nef and $A^2 \ne 0$. Then 
$$H^2(X,A+K_X) = 0.$$

\endclaim

\proof.  (A)
In a first step we show that we may assume $X$ to be $\bQ$-factorial.
(Actually, in our application in Section~5, it will be clear that we may
always assume $X$ to be $\bQ$-factorial, so the reader only interested
in the applications may skip (A)).
\smallskip
\noindent
In fact, if $X$ is not $\bQ$-factorial, there exists a bimeromorphic
map $f: Y \la X$ from a normal $\bQ$-factorial K\"ahler space with at
most terminal singularities ([Ka88,4.5$'$]). Moreover $f$ is an
isomorphism in codimension $1$ and $f$ is projective since $X$ has
only isolated singularities. Now consider the reflexive sheaf
$$\cH = f^*(\cO_X(A))^{**}.$$
Choose a number $r$ such that $A^{[r]}$ is locally free. Then 
$$ \cH^{[r]} = f^*(\cO_X(rA)), \eqno (1)$$
since both sheaves are reflexive and coincide in codimension $1$.
Thus $\cH$ is nef (as $\bQ$-line bundle) with $\cH^2 \ne 0$. 
Once we know the result in the $\bQ$-factorial case, we get
$$ H^2(Y,\cH \hotimes K_Y) = 0.$$
So by the Leray spectral sequence, we only have to show
$$R^1f_*(\cH \hotimes K_Y) = 0.$$
This however follows from [KMM87, 1-2-7]. Actually this citation deals
with the algebraic case. However first notice that our statement is
local around the isolated singularities of $X$.  Now isolated
singularities are algebraic by Artin's theorem, i.e.\ we can realize
an open neighborhood of an isolated singularity as an open set in a
normal algebraic variety. So locally on $X$ the map $f: Y \la X$
can be realized algebraically. Now we can approximate $\cH$ by
algebraic reflexive sheaves $\cH_k$ up to high order $k$ and then
apply [KMM87, 1-2-7] to get the vanishing $R^1f_*(\cH_k) = 0$. This
sheaf coincides with $R^1f_*(\cH)$ to high order, so $R^1f_*(\cH)$
vanishes to high order.  For $k$ approaching $\infty$, we obtain the
vanishing we are looking for.

\bigskip

\noindent
(B) From now on we assume $X$ to be $\bQ$-factorial.
We proceed as in the proof of 4.1. First of all choose $r$ such that 
$A^{[r]}$ is locally free. Then choose a singular metric $h$ with positive
curvature current on $A^{[r]}$. Now ${{1} \over {r}}h$ is a metric at least
on $A \vert X_{\rm reg}$ with positive curvature current $T$ extending to
all of $X$. We argue as in the first part of the proof of 4.1 to obtain
the divisor $D$ and the current $R$, however $D$ is only an integral Weil
divisor. By the same arguments as in 4.1 we can still reduce the problem
to proving
$$
H^2(D,\cO_D(A+K_X)) = 0.
$$
(Notice that ${\cal I}_D \otimes \cO_X(A+K_X) = \cO_X(-D+A+K_X)$
outside a finite set and that by definition $\cO_D(A+K_X) =
\cO_X(A+K_X) \vert D)$).  Now $D$ is Cohen-Macaulay; here we need in
an essential way that locally $X$ is the quotient of a hypersurface by
a finite group. To be more detailed, we can write locally $X = V/G$
with $V$ a hypersurface singularity and $G$ a finite group (see e.g.\ 
[Re87]). Let $\pi: V \la X$ be the quotient map and let $\hat D =
\pi^*(D)$.  If we can prove that $\hat D$ is Cohen-Macaulay, then $D$
will be Cohen-Macaulay, too, since this property is $G$-invariant. So
we may assume that $X = V$. Now $X$ is (locally) a compound Du Val
singularity [Re87], i.e.\ a $1$-parameter deformation of a
$2$-dimensional rational double point. Hence we can find a {\it
  Cartier} divisor $H \subset X$ through $x_0$ which has just a
rational double point at $x_0$. Now consider $D \cap H$. This is a
Weil $\bQ$-Cartier divisor on $H$. Since $x_0$ is a quotient
singularity of $H$, we can argue as above to see that $D \cap H$ is
Cohen-Macaulay. Hence $D$ has a hyperplane section through $x_0$ which
is Cohen-Macaulay. Thus $D$ is Cohen-Macaulay at $x_0$ itself.
\smallskip

\noindent Therefore we have by Serre duality 
$$
H^2(D,\cO_D(A+K_X)) \simeq \Hom(\cO_D(A+K_X),\cO_D(K_D)).
$$
Suppose that $H^2$ does not vanish. Then we obtain a non-zero
homomorphism $s: \cO_D(A+K_X) \la K_D$.  This $s$ must be generically
non-zero. In fact, $D$ is generically Gorenstein.  Hence
$\cO_D(K_D)_x$ is isomorphic to an ideal in $\cO_{X,x}$ for all $x$,
in particular $K_D$ has no torsion sections, $D$ being Cohen-Macaulay;
see [Ei95] in the algebraic case.  Let $X_0$ be the regular part of
$X$, this means that we eliminate a finite set from $X$, all
singularities being terminal.  Let denote $D_0 = D \cap X_0$ and let
$s_0 = s \vert X_0$.  Then by adjunction we have
$$ 0 \ne s_0 \in H^0(D_0,\cO_D(-A+D)).$$
From now on we argue as in 4.1 just working on $X_0$ instead of $X$.
The only exceptions are calculations of intersection numbers and
Hodge index arguments. Here we still need to argue on $X$ - we do not
have any problems with singularities since all divisors are $\bQ$-Cartier.
\qed

\section{\S5. The case $K_X^2 = 0, K_X \ne 0$}

The second ingredient for the proof of the abundance theorem for
K\"ahler threefolds is the following weak analogue of 4.3 in case
$K_X^2 = 0$ (however one should have in mind that we are dealing with
a cawe which does not exist a posteriori).

\claim 5.1\ Theorem| Let $X$ be a normal compact K\"ahler threefold with at
most terminal singularities such that $K_X$ is nef. Suppose that $K_X^2 = 0$, 
$K_X \ne 0 $, and that $X$ is simple and not Kummer. Then
$$ h^2(X,mK_X) \leq 1$$
for all $m \in \bN$.
\endclaim 

As already mentioned the essential property derived from $X$ being simple 
non-Kummer is that $\pi_1(X)$ is finite [Ca94].

\bigskip \noindent
{\bf 5.2 Start of the proof.} Using Kawamata's
$\bQ$-factorialisation theorem (compare with proof of 4.5), we may
assume that $X$ is $\bQ$-factorial.  Suppose $h^2(X,mK_X) \geq 2$.
Using Serre duality we get -- following Miyaoka and
Shepherd-Barron -- (many) non-split extensions
$$ 0 \la K_X \la \cE \la mK_X \la 0  \eqno (S)$$ 
with reflexive sheaves $\cE$ of rank $2$. 
We note
$$ c_1(\cE) = (m+1)K_X \eqno (5.2.1) $$
and
$$ c_2(\cE) = 0. \eqno (5.2.2)$$

\bigskip \noindent
{\bf 5.3 The unstable case.}\ 
\smallskip 

\noindent {\bf (5.3.a)} Here we will assume 
that every non-split extension $\cE$ as in (S) is $\omega $-
un\-sta\-ble for some fixed K\"ahler form $\omega$ independent of
$\cE$.  Let $\cA_S \subset \cE$ be the $\omega$-maximal destabilizing
subsheaf. Then $\cA_S$ is a $\bQ$-line bundle and we determine a
$\bQ$-line bundle $\cB_S$ such that $\cE/\cA_S = \cI_Z \cB_S$ with
some subspace $Z$ of codimension at least 2 (actually $Z$ is
generically (i.e.\ on the smooth part of $X$) locally a complete
intersection or finite and supported in ${\rm Sing}X)$.  Since $K_X
\ne 0$, we obtain injective maps
$$
\phi_S: \cA_S \la mK_X
$$
and
$$
\psi_S: K_X \la \cI \cB_S.
$$
Now there are (up to $\bC^*$) only finitely many maps $\phi_S:
\cA_S \la mK_X$ with some $\bQ$-line bundle $\cA_S$ arising as maximal
$\omega$-destabilizing subsheaf for some extension (S). In fact, fix
$\phi = \phi_S: \cA \la mK_X$.  Then by (6.13) there are only finitely
many maximal reflexive subsheaves $\cA' \subset mK_X$ such that $\cA'
\not \subset \cA$. So we may suppose $\cA' \subset \cA$. Then 
$$
\cA'\cdot \omega^2 < \cA \cdot \omega^2.
$$
Actually, putting
$$
\epsilon := {\rm min}Y_j \cdot \omega^2,
$$
where the minimum runs over the finitely many irreducible hypersurfaces
$Y_j \subset X$, we have
$$
\cA'\cdot \omega^2 \leq \cA \cdot \omega^2 - \epsilon.
$$
On the other hand, restricting ourselves to $\cA' $ of the form
$\cA' = \cA_{S'}$, we have by the destabilizing property
$$
\cA' \cdot \omega^2 \geq {{c_1(\cE) \cdot \omega^2} \over {2}}. \eqno (*)
$$
Having in mind that 
$$
\cA' = \cA \otimes \cO_X(-\sum \lambda_j Y_j),
$$
the finiteness of irreducible hypersurfaces in $X$ gives the
finiteness claim, since $(*)$ reads
$$
(\cA - \sum \lambda_jY_j) \cdot \omega^2 \geq {{c_1(\cE) \cdot
\omega^2} \over {2}}.
$$  
\smallskip 

\noindent So we have only finitely many possible maps 
$\phi$ (up to $\bC^*$).  In the same way (by dualizing) we have only
finitely many maps $\psi$ (up to $\bC^*$).  In (2) we prove that
$(\phi,\psi)$ and $(\lambda \phi,\lambda \psi)$ with $\lambda \in
{\Bbb C}^*$ always define isomorphic extensions (S).  Therefore in
total $(\lambda \phi,\mu \psi)$ with $\lambda, \mu \in {\Bbb C}^*$
define just a $1$-dimensional space of extensions, whence $h^2(X,mK_X)
\leq 1$.  \smallskip

\noindent 
{\bf (5.3.b)} We shall now prove that the extension class defining (S)
is already determined by $\phi$ and $\psi$ (modulo $\bC^*$).  So take
another extension
$$
0 \la K_X \la \cE' \la mK_X \la 0 \eqno (S')
$$
with the same destabilizing sheaves $\cA$ and $\cB$ and with the
same morphisms $\phi$ and $\psi$ (the case of $(\phi,\psi), (\lambda
\phi, \lambda \psi)$ is exactly the same).  Let $D$ be the divisorial part of
$$
\{\phi = 0\} \cup \{\psi = 0\}\cup {\rm Sing}(X);
$$
then we obtain a splitting of the sequence $(S)$  over $X \setminus D$ 
via $\phi:$
$$
\cE \simeq mK_X \oplus K_X \simeq \cA \oplus \cB,
$$
and an analogous splitting of $\cE'$ over $X \setminus D$. Observe also 
that $\cA = mK_X - D$ and $\cB = K_X + D$. Thus we obtain an isomorphism 
$$ f: \cE \la \cE'$$
on $X \setminus D$ making the two extensions $(S)$ and $(S')$ isomorphic
over $X \setminus D:$
$$ \matrix{ 0 & \la & K_X & \la & \cE & \la & mK_X & \la & 0 \cr
              & &  \Arrowvert &    & \downarrow  &  & \Arrowvert & &  \cr
            0 & \la & K_X & \la & \cE' & \la & mK_X & \la & 0. \cr } $$
          
It remains to extend the map $f$ to $X$.  Let us notice that we may
assume $Z = \emptyset$. In fact let $Z_1$ be the codimension 1 part of
$Z$. Restricting our two exact sequences describing $\cE$ to $D$, we
see that (modulo torsion at finitely many points)
$$
\cA \vert D = K_X \vert D,
$$
hence
$$
(m-1)K_X \cdot D = D \cdot D. \eqno (5.3.1)
$$
In particular we note that $D \vert D$ is nef, hence $D$ itself is nef.
Now (5.3.1) yields
$$
0 = c_2(\cE) = c_1(\cA) \cdot c_1(\cB) + c_2(\cI_Z B) = (mK_X-D)
\cdot (K_X+D) + Z_1 = Z_1,
$$
hence $Z_1 = \emptyset$. In particular $Z \subset \Sing X$ has
codimension at least 3.  \smallskip \noindent This shows that we may
ignore $Z$ in all our following considerations; in what follows
restriction will always that we also divide by torsion.  \smallskip
\noindent Take a local section $s \in \cE(U)$ over a small disc $U$.
We need to show that $f(s) \in \cE'(U);$ a priori we only know $f(s)
\in \cE'(D)(U)$.  Let
$$
\kappa : \cE \to mK_X, \ \kappa' : \cE' \to mK_X
$$
and
$$
\lambda: \cE \to \cB, \ \lambda' : \cE' \to \cB
$$
be the canonical maps. 

Now consider the exact sequence
$$
0 \la \cE \ {\buildrel {\kappa \oplus \lambda} \over {\la}} \ mK_X
\oplus \cB \ {\buildrel {\rho} \over {\la}} \ mK_X \vert D \la 0
$$
Here $\rho (u \oplus v) = u_D - \tau(v_D)$, where $\tau: \cB_D \to mK_X \vert D$ is the canonical sequence arising by
restricting our sequences and the maps $\phi$ and $\psi$ to $D$. 
Analogously for $\cE'$. Suppose we know
$$
h^0(D_c,\cO_{D_c}) = 1 \leqno (+)
$$
for any connected component $D_c$ of $D$. Actually it suffices to
know this for $D_c \cap {\rm reg}X$.  Then $\tau'^{-1} \circ \tau = a
\ {\rm id }, $ and of course we can normalize (in the extension class)
to $a = 1$ (our arguments are local around every individual connected
component of $D$).  Therefore we can construct a diagram
$$
\matrix{ 0 & \la & \cE & \la & mK_X \oplus \cB & \to & mK_X \vert D
  & \la & 0 \cr & & \downarrow & & \Arrowvert & & \downarrow & & \cr 0
  & \la & \cE' & \la & mK_X \oplus \cB & \la & mK_X \vert D & \la & 0.
  \cr }
$$
Here the map $f: \cE \to \cE'$ is defined only on $X
\setminus D$ and meromorphic on $X$, the left square is commutative on
$X \setminus D$ and the right square is commutative on $X$. It is
immediately checked that $mK_X \vert D \to mK_X \vert D$ is the
identity (consider images of elements $u \oplus 0$), hence $f(s) \in
\cE'(U) $ and thus $f$ extends to a global isomorphism making the two
extensions isomorphic. \medskip \noindent

\noindent {\bf (5.3.c)} It remains to prove (+) and we may assume $D = D_c$ for simplicity of notations. So suppose that
$h^0(D,\cO_D) \geq 2$, resp.\ $h^0(D\cap{\rm reg}X, cO_{D \cap 
{\rm reg}X}) \geq 2$.  Since $H^1(X,\cO_X) = 0$, we obtain
$$
H^1(X,\cO_X(-D)) \ne 0.
$$
Let $X_0 $ be the regular part of $X$ and $D_0 = D \vert X_0$.
Since the singularities of $X$ are at most finite,
$$
H^1(X_0,\cO(-D_0)) = H^1(X,\cO(-D)) \ne 0.
$$
Hence ${\rm Ext}^1(\cO_{X_0},\cO(-D_0)) \ne 0$, and therefore there
exists a non-split extension
$$
0 \la \cO_{X_0} \la \cF_0 \la \cO(D_0) \la 0, \leqno (E_0)
$$
with a locally free sheaf $\cF_0$ over $X_0$.  Now $\cF_0$ has a
unique reflexive extension to $X:$ consider a singular point $x_0 \in
X$ and let $U$ be a Stein neighborhood of $x_0$. Then
$$
H^1(U \setminus x_0,\cO(-D_0)) = H^1(U,\cO(-D)) = 0,
$$
hence $(E_0)$ splits over $U \setminus x_0:$
$$
\cF_0 \simeq \cO \oplus \cO_{U \setminus x_0}(D_0).
$$
Hence $\cF_0$ extends to a reflexive sheaf $\cF$. Moreover $(E_0)$ extends to 
$$
0 \la \cO_X \la \cF \la \cO_X(D) \la 0. \leqno (E)
$$
In particular ${\rm Ext}^1(\cO_X(D),\cO_X) \ne 0$. This is easily
seen to be equivalent to
$$
{\rm Ext}^1(K_X + \cO_X(D),\cO_X(D)) \ne 0,
$$
hence $H^2(X,K_X+ \cO_X(D)) \ne 0$ by Serre duality. Thus 
$$
D^2 = 0
$$
by (4.3).\smallskip

\noindent
We observe $c_1(\cF)^2 = c_2(\cF) = 0$, therefore $\cF$ cannot be
$\omega$-stable. So let $\cA'$ be the maximal $\omega$-destabilizing
subsheaf, and we obtain as before for $\cE$ a sequence
$$
0 \la \cA' \la \cF \la \cI_{Z'}\cB' \la 0, \leqno (F)
$$
where $Z'$ has generically dimension 1 or is contained in the
singular locus of $X$. As before, $\cA' \subset \cO_X(D)$ and $\cO_X
\subset \cB'$ so that there is an effective divisor $D'$ such that
$\cB' = \cO_X(D') $ and $\cA' = \cO_X(D-D')$. By $(E)$ and the fact
that $a(X) = 0$ and $H^1(X,\cO_X) = 0$, we deduce that $h^0(X,\cF) = 2$. 
Hence $(F)$ yields
$$
h^0(X,\cA') = 1.
$$
So we can write $\cA' = \cO_X(E)$ with an effective divisor $E$.
Thus $D = D'+E$. By restricting $(E)$ and $(F)$ to $D'$, we obtain
$\cA' \vert D' = E \vert D' = \cO_{D'}$ and that the 1-dimensional
part $Z_1'$ of $Z'$ is empty, so that $Z'$ is at most finite.
Alternatively, we calculate
$$
0 = c_2(\cF) = E \cdot D' + Z_1'
$$
and conclude together with the observation that $D' \vert D'$ and
thus $D'$ is nef.
\smallskip 

\noindent So we have a decomposition $D=D'+E $ with $D$ and $D'$ nef. 
By instability, $E \ne 0$.  Having in mind that $D$ is connected, we
are going to prove that $D$,$D'$ and $E$ are proportional (even
numerical proportionality would be sufficient, which in our situation
($X$ simply connected with $a(X) = 0$) gives equality).  \smallskip
\noindent Assuming this proportionality for the moment, we obtain $D'
= aE$ and $D = (a+1)E$.  Since $\cA'$ destabilizes, we have $a \leq 1$.
By restricting the sequence $(E)$ and $(F)$ to $D'$ we obtain:
$$
\cO_{D'}(E) = \cO_{D'} \eqno (*)
$$
up to torsion. For the simplicity of notation we suppress the
torsion and agree, when taking a restriction, that we also divide by 
the torsion. We
may assume that $h^0(D',\cO_{D'}) = 1$. In fact, if $h^0(D',\cO_{D'}) \geq
2$, then we substitute $D$ by $D'$ and argue as before; of course this
procedure has to terminate after finitely many steps. So
$h^0(D',\cO_{D'}) = 1$ and consequently
$$
h^0(D',N_{D'/D}^{*\mu}) \ne 0
$$
for some positive integer; hence $h^0(D',N_{D'/X}^{*\mu}) \ne 0$.
(We may neglect the torsion, because we may compute over $X_0$).
\smallskip \noindent Now this non-vanishing implies that $\mu D' \vert
D_0$ is trivial, where $D_0$ is the reduction of $D'$. In fact, let
$\mu D' = \sum a_i Y_i$ and take a non-zero section $s$ of $\cO_{\mu
  D'}(\mu D')$. Let $ \tilde D = \sum b_i Y_i \subset \mu D'$ be the
maximal subdivisor such that $s \vert \tilde D = 0$.  Introducing 
$c_i = a_i - b_i \geq 0$, we obtain a section 
$$
s' \in H^0(D_0,\cO_{D_0}(-\sum c_i Y_i))
$$
such that $s' \vert Y_i \ne 0$ for
all $i$. Fix a K\"ahler form $\omega$ and let $\alpha_{ij} = Y_i \cdot
Y_j \cdot \omega $ for $i \ne j$. Then we obtain for all $j$
$$
- \sum_i c_i Y_i \cdot Y_j \cdot \omega = - \sum c_i \alpha_{ij} \geq 0.
$$
Since $D^2 = 0$ and $D$ is nef, we have $D \cdot Y_j = 0$ for all $j$ and 
therefore 
$$
Y_j^2  = - {{1} \over {a_j}} \sum_{i \ne j} a_i Y_i \cdot Y_j
$$
so that we arrive at the inequalities (for each $j$)
$$
\sum_{i \ne j} a_ic_j \alpha_{ij} \geq \sum_{i \ne j} a_j c_i \alpha_{ij}.
$$
By simple algebraic considerations this is only possible if we have
always equality. This means that the divisor $D^* = \sum c_i Y_i$
fulfills $D^* \cdot Y_j = 0$ for all $j$. Hence $D^*$ is nef, and the
proportionality arguments below shows that $D^* = cD'$ for some
positive number $c$.  \smallskip \noindent Because of the
non-vanishing of $s'$ on $Y_j$, $D^* \vert Y_j $ is trivial, hence
$D^* \vert D_0$ is trivial. Suppose for the moment that $\mu = 1$.
Then $D^* \subset D'$ so that $c < 1$ and in total we have $D^* \vert
D_0$ trivial, and $D^* \vert D^* $ torsion by $(*)$. Hence [Mi88a, 4.1]
says that $D^* \vert D^* $ is trivial.  Now the exact sequence
$$
0 \la \cO_X \la \cO_X(D^*) \la \cO_X(D^*)\vert D^*  \la 0
$$
implies by $(*)$ -- keeping in mind $H^1(X,\cO_X) = 0$ -- that
$h^0(X,\cO_X(D')) \geq 2$, contradicting $a(X) = 0$.  \smallskip
\noindent So we are left with the case $\mu \geq 2$. We deal with $\mu
= 2$ and leave the trivial modifications in the general case to the
reader. The difficulty here is that possibly $c > 1$ so that $D'
\subset D^*$, otherwise we conclude as before.  At least we know that
$D^* \subset \mu D'$ and we are going to show that $\mu E \vert \mu D'
$ is trivial; then we are done.  This does not follows directly from
restricting $(E)$ and $(F)$; instead we take $S^2$ and obtain an
injection $\cO_X(2E) \la S^2(\cF)$. Restricting to $2D'$, we obtain a
non-zero map $\cO_{2D'}(2E) \to \cF \vert 2D'$.  Then either the
induced map $\cO_{2D'}(2E) \to \cO_D(D)) $ is non-zero; this implies
$H^0(D',\cO_{D'}(-D')) \ne 0 $ so that $h^0(X,\cO_{2D'}) > 1$ and we may
take $D = 2D'$ whence $\mu = 1$. Or this map vanishes; then we get a
non-zero map $ \cO_{2D'}(2E) \to \cO_{2D'}$. This map is an
isomorphism and settles our claim.

\medskip \noindent {\bf (5.3.d)} It remains to settle the
proportionality.  After passing to a desingularization this comes down
to prove the following statement. Let $X$ be compact K\"ahler, $A =
\sum_{i=1}^s a_i Y_i$ and $B = \sum_{i=1}^s b_i Y_i$ be effective
divisors, $a_i$ and $b_i$ positive, with connected supports. Suppose
that $A$ and $B$ are nef, and that $A \cdot Y_i = B \cdot Y_i = 0$ for
all $i$ (in particular $A^2 = B^2 = A \cdot B = 0$). Then $A = cB$
(not only for numerical equivalence).  Observe that if $X$ is a
surface, then this is nothing than Zariski's lemma, which is usually
formulated for fibers of maps to curves, but which works in this
context; therefore the claim also follows for projective manifolds by
taking hyperplane sections and applying Lefschetz. If $X$ is merely
K\"ahler, then we consider the vector space $V \subset H^2(X,\bR)$
generated by the classes of the hypersurfaces $Y_i \subset X$. Let $W$
be the direct sum $\bigoplus \bR \cdot Y_i;$ and let $Q$ be the
bilinear form
$$Q(Y_i,Y_j) = -Y_i \cdot Y_j  \cdot \omega.$$
In this situation we apply [BPV84,lemma I.2.10] to conclude.  
\qed

\bigskip \noindent
{\bf 5.4 The stable case.} \ By 5.3 we are reduced to the case that 
some extension $\cE$ is $\omega$-stable for some K\"ahler form 
$\omega$. By (5.2.1) and (5.2.2) we have in particular
$$
c_1^2(\cE) \cdot \omega = 4 c_2(\cE) \cdot \omega,$$
which implies
that $\cE$ is projectively flat, at least on the regular part $X_0$.
In fact, this is well-known if $X$ is smooth and $\cE$ is locally
free. But the proof generalizes to our case since the singularities of
$X$ and $\cE$ are in codimension at least 3. Now we follow the
arguments in [Ko92,p.113/114].  \smallskip

\noindent Assume first that the degree of finite \'etale covers of 
$X_0$ is bounded: $\pi_1^{\rm alg}(X_0) $ is finite. After performing a 
finite \'etale cover, we may assume $\pi_1^{\rm alg}(X_0) = 0 $.  Since 
$\cE\vert X_0$ is projectively flat, $\cE^* \otimes \cE \vert X_0$ is
hermitian flat and therefore given by a unitary representation $\rho$
of $\pi_1(X_0)$.  Since $\rho(\pi_1(X_0))$ is residually finite, it
follows that $\rho$ is trivial, hence $\cE^* \otimes \cE$ is trivial.
This implies, using the exact sequence
$$
0 \la K_X \la \cE \la mK_X \la 0
$$
and dualizing that (over $X_0)$ $h^0( \cE^* \otimes K_X) = 0$, that
$h^0(\cE^* \otimes mK_X) \geq 4$ and that therefore $h^0((m-1)K_X)
\geq 3$.  This contradicts $a(X) = 0$.  \smallskip

\noindent If  $\pi_1^{\rm alg}(X_0) $ is infinite, we just take over 
the arguments of [Ko92p.114]: since the local fundamental groups of $X$ at
the singularities are finite, any finite \'etale cover $h$ of $X_0$ of
sufficiently large degree extends to a covering $h: \tilde X \to X$
which can be written in the form $h = g \circ f$, where $f: \tilde X
\to X'$ is \'etale and $g: X' \to X$ is \'etale outside the singular
locus.  Therefore $\pi_1(X')$ is infinite, contradicting the fact that
$X'$ is simple non-Kummer.  \qed
\bigskip

\section{\S6. The inequality $K_X \cdot c_2(X) \geq 0$}

\vskip10pt

The aim of this section is to prove

\claim 6.1\ Theorem| Let $X$ be a minimal K\"ahler 3-fold with $a(X) = 0$. 
Then $K_X \cdot c_2(X) \geq 0$. \endclaim

This inequality is an important step in the proof of abundance for
K\"ahler threefolds.  In the projective case, it follows from
Miyaoka's inequality $K_X^2 \leq 3c_2(X)$ which in turn is a
consequence of his generic nefness theorem for the cotangent bundle
(relying on char $p$ methods). \smallskip

\noindent
The rest of this section consists of the proof of 6.1 together with
some auxiliary propositions (6.9, 6.10, 6.12/6.13).

\bigskip \noindent {\bf 6.2}\ {\bf Reduction to the unstable case.}
Suppose that there is a sequence $(\omega_j)$ of K\"ahler metrics
converging in $H^2(X,\Bbb R)$ to $K_X$ such that $T_X$ is
$\omega_j$-stable for all $j$.  Then we have
$$
c_1^2(X) \cdot \omega_j \leq 3c_2(X) \cdot \omega_j
$$
for all $j$ by Proposition 6.9. Taking limits, we obtain
$$
K_X^3 \leq 3 K_X \cdot c_2(X),
$$
hence our claim results from $K_X^3 = 0$. 
\smallskip

\noindent
So from now on we shall assume that {\it $T_X$ is $\omega$-unstable
for all $\omega$ near $K_X$ (in $H^2(X,\Bbb R)).$}

\bigskip \noindent
{\bf 6.3} {\bf The setup}
Let $\cS_{\omega}\subset T_X$ be the maximal destabilizing subsheaf with 
respect to $\omega$. Let $r$ denote its rank. Then by Corollary 6.13 below, 
there are only finitely many choices for $\cS_{\omega}$, hence there 
exists an open set in the
K\"ahler cone of $X$ having $K_X$ as boundary point such that
$\cS_{\omega}$ does not depend on $[\omega]$ for $[\omega] \in U$. 
We shall write
$$ S = \det \cS$$
and let $Q = T_X/\cS$, a torsion free sheaf of rank $1$ or $2$.
We notice
$$ c_1(Q) = -K_X - S  \eqno (6.3.1)$$
and, if $ r = 1$, 
$$ c_2(Q) = c_2(X) + S \cdot (K_X+S). \eqno (6.3.2)$$
$c_3(Q)$ will be irrelevant for us.  
The instability of $T_X$ gives
$$
S \cdot \omega^2 \geq {{-K_X \cdot \omega^2} \over {3}}, \eqno (6.3.3)
$$
for $\omega \in U$ in case $ r = 1$ and
$$
{{S \cdot \omega^2} \over {2}} \geq {{-K_X \cdot \omega^2} \over
{3}} \eqno (6.3.3a)
$$
in case $r = 2$. 
We claim 
$$
K_X^2 \cdot S = 0. \eqno (6.3.4)
$$
In fact, (6.3.3) resp.\ (6.3.3a) gives in the limit $K_X^2 \cdot S
\geq 0$. Since we may assume $K_X^2 \ne 0$, the tangent sheaf $T_X$ is
$K_X$-semi-stable by Enoki [En87]. This implies $K_X^2 \cdot S \leq 0$, 
hence (6.3.4) follows.
\smallskip

The next lemma is a general statement on K\"ahler 3-folds with $a(X) = 0$,
independent from our setup.

\claim 6.4\ Lemma| Let $X$ be a normal compact K\"ahler 3-fold.  Let
$A$ and $B$ be $\bQ$-line bundles on $X$ and let $A$ be nef with $A^2
\ne 0$.  If $A^2 \cdot B = 0$, then $A\cdot B^2 \leq 0$.  \endclaim

\proof. By passing to a desingularization we may assume $X$ smooth.
Fix a K\"ahler class $\omega $ and apply (2.5) with $\alpha = c_1(B),
\beta = K_X + \varepsilon \omega $ and $\gamma = K_X$.  Then expand
in terms of powers of $\varepsilon$ to obtain the claim.\qed
 
\claim 6.4.a\ Corollary| In our setup (6.3) we have $K_X \cdot S^2
\leq 0$ if $K_X^2 = 0$.  \endclaim

\proof. This follows from 6.4 via 6.3.4  \qed

Lemma 6.3 is of course not true in case $A^2 = 0$. Thus in order to
obtain (6.4.a) also in case $K_X^2 = 0$ we need more specific
arguments:

\claim 6.5\ Lemma| Let $X$ be a simply connected minimal K\"ahler
3-fold with $a(X) = 0$ and $K_X^2 = 0$.  Let $L$ be a $\bQ$-line
bundle on $X$.  Then $K_X \cdot L^2 \leq 0$.
\endclaim

\proof. Assume that $K_X \cdot L^2 > 0$. If a positive integer $c$
satisfies the following condition:
$$
2c^2 K_X \cdot L^2 > -K_X \cdot c_2(X) \eqno (*),
$$
then by Riemann-Roch we easily get asymptotically
$$
\chi(X,mK_X+cL) \sim m.
$$
Observe also that $(*)$ is satisfied for large $c$ since $K_X \cdot
L^2 > 0$ by assumption.  So let us fix such a number $c$. Then we
conclude
$$
h^2(mK_X+cL) \geq C m. \eqno (**)
$$
In fact, otherwise $h^0(mK_X+cL) \geq C m$ by $(*)$, contradicting $a(X) = 0$. 
\smallskip 

\noindent
Now, as in section 5, we obtain ``many'' extensions
$$
0 \la K_X+cL \la \cE \la mK_X \to 0.
$$
Observe that $\cE$ cannot be $\omega$-stable for $\omega$ near
$K_X$. In fact, in that case we had
$$
c_1(\cE)^2 \cdot \omega \leq 4 c_2(\cE) \cdot \omega,
$$
hence $c_1(\cE)^2 \cdot K_X \leq 4 c_2(\cE) \cdot K_X$ in the
limit. This comes down to
$$
c^2 K_X \cdot L^2 \leq 0,$$ contradicting our assumption.
\smallskip

\noindent 
We proceed exactly in the same way as in section~5, introducing the
divisors $D_m$, and now $(**)$ and the arguments in section~5 yield
$$
h^0(\cO_{D_m}) \geq Cm,
$$
for large $m$. On the other hand, 
$$
(m-1)K_X + cL^* \vert D_m = D_m \vert D_m,
$$
again referring to section~5, hence for large $m$, the normal
bundle $N_{D_m}$ gets more and more ``nef''. However, to have many
functions on $D_m$ means to have a tendency to negativity for the
normal bundle.  So we will show that (+) and (++) are contradictory.
By passing to a subsequence - having in mind that $X$ carries only
finitely many irreducible hypersurfaces - we can suppose the
following.
$$
D_m = \sum_{i=1}^s a_{m,i}Y_i  + \sum_{s+1}^t a_j Y_j,
$$
where $a_{m,i} < a_{m+1,i}$ and the $a_j$ are independent of $m$.
Put $R = \sum_{s+1}^t a_j Y_j$ and $Y = \sum_1^s Y_i$.  Then
$$
N^*_{D_m/D_{m+1}} = (\sum_{i=1}^s a_{m,i}Y_i + R) \vert Y.
$$
Since by (+), $h^0(N^*_{D_m/D_{m+1}}) > 0$, the sequence of
divisors $ - \sum a_{m,i} Y_i \vert Y$, suitably normalized, converges
to an effective non-zero divisor on $Y$.  Thus $N^*_{D_m\vert X} \vert Y$, 
suitably normalized converges to an effective non-zero divisor on
$Y$.  On the other hand, its dual is nef by (++). This is a
contradiction.  \qed

\claim 6.5.a\ Corollary|  In our setup (6.3) we have $K_X \cdot S^2 \leq 0$.
\endclaim 

\noindent
{\bf 6.6 The Case:} \ $\rk \cS = 1$ {\bf and} $Q$ \ {\bf stable.}\ 
By ``$Q$ stable'' we mean that there is a sequence of K\"ahler forms
$(\omega_j)$ converging to $K_X$ (as classes) such that $Q$ is
$\omega_j$-stable for all $j$.  Then by Proposition 6.9 we have
$$ c_1^2(Q) \cdot \omega_j \leq 4c_2(Q) \cdot \omega_j,$$
hence
$$ c_1^2(Q) \cdot K_X \leq 4c_2(Q) K_X.$$
Putting in (6.3.1) and (6.3.2) we obtain
$$ K_X \cdot (K_X+S)^2 \leq (4c_2(X) + 4S \cdot (K_X+S)) \cdot K_X,$$
which in turn yields
$$ K_X \cdot c_2(X) \geq - {{3}\over {4}}K_X \cdot S^2.$$
Thus 6.4.a gives $K_X \cdot c_2(X) \geq 0$.

\bigskip \noindent {\bf 6.7 The Case:}\ $\rk \cS = 1$ {\bf and} $Q$
{\bf is unstable}.\ 
After the previous case it is clear what unstable
has to mean: $Q$ is $\omega$-unstable for all $\omega$ near $K_X$.
Then we obtain a destabilizing sequence
$$
0 \la L_1 \la Q \la \cI_B L_2 \la 0
$$
where $L_i$ are reflexive of rank $1$ and $\dim B \leq 1$.  This
sequence is - as usual - independent of $\omega$, if $\omega$ is
sufficiently near to $K_X$ and contained in a suitable open set $U$ as
in (6.3).  We first claim
$$
K_X^2 \cdot L_1\leq 0. \eqno (6.7.1)
$$
To verify this, let $\cR$ be the cokernel of
$$
T_X \la \cI_B L_2 \la 0.
$$ 
Then we have an exact sequence 
$$
0 \la S \la \cR \la L_1 \la 0.
$$
Of course we may assume $K_X^2 \ne 0$. Then by Enoki [En87], $T_X$ is 
$K_X$-semi-stable, hence 
$$
c_1(\cR) \cdot K_X^2 \leq 0.
$$
This implies $K_X^2 \cdot (L_1+S) \leq 0$ by the last exact
sequence. Now (6.3.4) gives our claim (6.7.1) \medskip 

\noindent Next we show
$$
K_X \cdot L_1 = 0. \eqno (6.7.2)
$$
In fact, the destabilizing property for $L_1$ reads
$$
L_1 \cdot \omega^2 \geq {{c_1(Q) \cdot \omega^2} \over {2}},
$$
hence
$$
L_1 \cdot K_X^2 \geq {{c_1(Q) \cdot K_X^2} \over {2}} = {{1} \over {2}}
(-K_X-S) \cdot K_X^2 = 0.
$$
We now conclude by (6.7.1).

\medskip \noindent Since $c_1(Q) \cdot K_X^2 = 0$, we also have 
$$
K_X^2 \cdot L_2  = 0. \eqno (6.7.3)
$$
Thus Lemma 6.4 applies:
$$
K_X \cdot L_i^2 \leq 0 \eqno (6.7.4)
$$
for $i=1,2$.\smallskip 

\noindent The final preparation is
$$
K_X\cdot L_1 \cdot L_2 = {{1} \over {2}}(K_X \cdot S^2 - K_X \cdot
L_1^2 - K_X \cdot L_2^2).
$$
This follows from the two equations
$$
K_X \cdot c_1^2(Q) = K_X \cdot (L_1+L_2)^2
$$
and 
$$
K_X \cdot c_1^2(Q) = K_X \cdot (K_X+S)^2 = K_X \cdot S^2.
$$
After all these preparations we conclude using (6.7.5) as follows.
$$
\eqalign{
K_X \cdot c_2(X)
&= K_X \cdot c_2(Q) + K_X \cdot S \cdot c_1(Q)\cr
&= K_X \cdot c_2(\cI_B) + K_X \cdot L_1 \cdot L_2 -K_X \cdot S^2\cr
&= K_X \cdot c_2(\cI_B) + {{K_X \cdot S^2} \over {2}}
- {{K_X \cdot L_1^2} \over {2}} - {{K_X \cdot L_2^2} 
\over {2}} - K_X \cdot S^2.\cr}
$$
Since $K_X \cdot c_2(\cI_B) \geq 0$ by nefness of $K_X$ we conclude
by virtue of (6.3) and (6.7.4).
\bigskip 

\noindent
{\bf 6.8 The Case:} \ $\rk \cS = 2$. \smallskip

\noindent
In this case we consider the maximal destabilizing subsheaf $Q^*
\subset \Omega^1_X = (T_X)^*$. Here it is convenient to switch
completely the notations: we denote the maximal destabilizing subsheaf
of $\Omega^1_X$ again by $\cS$ and let $Q$ denote the quotient. Then
$$
c_1(Q) = K_X - S
$$
and 
$$
c_2(Q) = c_2(X) - S\cdot (K_X-S).
$$
Now (6.3) yields $K_X^2 \cdot S = 0$. Applying again 6.4 gives $K_X
\cdot S^2 \leq 0$.  Now (6.6) and (6.7) run in completely the same
way; notice that some minus signs are irrelevant because $K_X^2 \cdot
S = 0$.

\claim 6.9\ Proposition| Let $X$ be a normal compact K\"ahler $n$-fold with
$\codim {\rm Sing}(X) \geq 3$. Suppose $a(X) = 0$.
Let $\omega$ be a K\"ahler form on $X$ and $\cE$ a torsion free coherent 
sheaf on $X$ of rank $r \geq 2$. If $\cE$ is $\omega$-stable, then
$$
c_1^2(\cE) \cdot \omega^{n-2} \leq {{2r} \over {r-1}}  c_2(\cE) \cdot
\omega^{n-2}.
$$
\endclaim

\proof. For simplicity of notations set $\mu = {{2r}\over {r-1}}$. 
\smallskip 

\noindent
(1) First we reduce the problem to the case ``$\cE$ reflexive''. So
suppose we know the assertion for reflexive sheaves and let $\cE$ be
torsion free.  Then we consider the quotient sheaf
$$
Q = \cE^{**}/\cE,
$$
which is supported on a complex subspace $Z \subset X$ of
codimension at least $2$.  Now
$$
c_2(Q) = -m c_2(\cI_Z),
$$
for some positive $m$, and $c_2(\cI_Z)$ is an effective cycle supported 
on $Z$, hence
$$
\omega^{n-2} \cdot c_2(Q) \leq 0. \eqno (*)
$$
Now $c_2(\cE^{**}) = c_2(\cE) + c_2(Q)$, hence $(*)$ implies
$$
c_2(\cE) \omega^{n-2} \geq c_2(\cE^{**}) \cdot \omega^{n-2}. \eqno (**)
$$
Notice that $\cE^{**}$ is stable because $\cE$ is ([Ko87,V.7.7]),
hence by our assumption
$$
c_1^2(\cE^{**}) \omega^{n-2} \leq \mu c_2(\cE^{**}) \cdot \omega^{n-2}.
$$
Since $c_1(\cE^{**}) = c_1(\cE)$, the inequality $(**)$ implies our
claim follows.\medskip 

\noindent 
(2) From now on we shall assume $\cE$ reflexive.  Choose a
desingularization $\pi: \hat X \la X$ by a sequence of blow-ups whose
centers all ly over the singularities of $X$ and $\cE$. Moreover we
may assume that $\hat \cE = \pi^*(\cE)^{**}$ is locally free (see
[GR71]). Let $\hat \omega = \pi^*(\omega)$. By definition of K\"ahler
forms on singular spaces $\hat \omega$ - which a priori exists only on
a Zariski open part of $\hat X$ - extends to a semipositive
$(1,1)$-form on all of $\hat X$.  We claim
$$
\cE {\rm \ is \ } \hat \omega-{\rm stable}. \eqno (+)
$$
Indeed, assume we have a subsheaf $\hat \cS \subset \hat \cE$ of
rank $s$ with
$$
{{c_1(\hat \cS) \cdot \hat \omega^{n-1}} \over {s}} \geq {{c_1(\hat \cE) 
\cdot \hat \omega^{n-1}} \over {r}}.
$$ 
Then consider 
$$
\cS = \pi_*(\hat \cS) \subset \pi_*(\hat \cE).
$$
Since $\pi_*(\hat \cE) $ is torsion free and since $\cE$ is
reflexive, we have $\pi_*(\hat \cE) \subset \cE$, hence $\cS \subset
\cE$.  Now
$$
c_1(\hat S) \cdot \hat \omega^{n-1} = c_1(\pi_*(\hat \cS)) \cdot 
\omega^{n-1} = c_1(\cS) \cdot \omega^{n-1},
$$
and
$$
c_1(\hat \cE) \cdot \hat \omega^{n-1} = c_1(\cE) \cdot \omega^{n-1},
$$
hence
$$
{{c_1(\cS) \cdot \omega^{n-1}} \over {s}} \geq {{c_1(\cE) \cdot 
\omega^{n-1}} \over {r}},
$$
contradicting the $\omega$-stability of $\cE$. This proves (+). 
\smallskip 

\noindent
Now $\hat \omega$ has the disadvantage not to be a K\"ahler form, but
it is on the boundary of the K\"ahler cone. To circumvent this
difficulty, let $E_i$ denote the exceptional components of the
exceptional set of $\pi$, then we can chose $a_i < 0$, such that 
$E := \sum a_i E_i$ is $\pi$-ample.  Thus
$$
\hat \omega_{\epsilon} := \hat \omega + \epsilon E
$$
is a K\"ahler class for all small positive $\epsilon$.  We claim
that $\hat \cE$ is $\hat \omega_{\epsilon}$-stable for $\epsilon$
small enough.  Indeed, suppose the contrary. Then there exists a
sequence $\epsilon-k$ converging to $0$ such that $\hat \cE$ is not $
\hat \omega_{\epsilon_k}$-stable. Let $\cS_i \subset \hat \cE$ be the
maximal destabilizing subsheaf with respect to $ \hat
\omega_{\epsilon_i}$. Since $a(\hat X) = 0$, we find $i_0$ such that
$\cS_i = \cS_j $ for all $i,j \geq i_0$ (Prop. 6.12), possibly after
passing to a subsequence (but even this could be avoided). So let 
$\cS=\cS_i$, $i \geq i_0$.  Then we have
$$
{{c_1(\cS) \cdot \hat \omega_{\epsilon_k}^{n-1}} \over {s}} \geq
{{c_1(\hat \cE) \cdot \hat \omega_{\epsilon_k}^{n-1}} \over {r}},
$$
so passing to the limit,
$$
{{c_1(\cS) \cdot \hat \omega^{n-1}} \over {s}} \geq {{c_1(\hat \cE)
\cdot \omega^{n-1}} \over {r}}.
$$
This contradicts (+).\smallskip 

\noindent Thus $\hat \cE$ is $\hat \omega$-stable for small
positive $\epsilon$.  Therefore $\hat \cE$ is Hermite-Einstein with
respect to $ \hat \omega_{\epsilon}$ and hence
$$
c_1^2(\hat \cE) \cdot \hat \omega_{\epsilon}^{n-2} \leq \mu \ 
c_2(\hat \cE) \cdot \hat \omega_{\epsilon},
$$
hence
$$
c_1^2(\hat \cE) \cdot \hat \omega^{n-2} \leq
\mu \  c_2(\hat \cE) \cdot \hat \omega^{n-2}.
$$ 
Since $\codim ({\rm Sing}(X) \cup {\rm Sing}(\cE)) \geq 3$, we conclude
$$
c_1^2(\cE) \omega^{n-2} \leq \mu \  c_2(\cE) \cdot \omega^{n-2}.
$$\qed

\claim 6.10\ Proposition| Let $L$ be a line bundle or a reflexive sheaf
of rank $1$ on the normal compact complex space $X$. Suppose $a(X) = 0$.
Let $\cS_i \subset L$
be reflexive subsheaves, $i \in I$. Then for all $i$ there are only finitely
many $j$ such that $\cS_j \not \subset \cS_i$.
\endclaim

\proof.  Of course we may assume $X$ smooth. Since $a(X)= 0$, the
complex space $X$ has only finite many irreducible hypersurfaces $Y_1,
\ldots Y_r$, therefore we can write
$$
\cS_i = L - \sum_{j=1}^r a_j^{(i)}Y_j
$$
with $a_j^{(i)} \geq 0$. Thus the claim is clear.  
\qed

\claim 6.11\ Definition| Let $\cF$ be a torsion free coherent sheaf on
a normal compact complex space and let $\cS \subset \cF$ be a
reflexive subsheaf with $0 < {\rm rk}\cS < {\rm rk}\cF$. We say that
$\cS$ is maximal, if there is no proper reflexive subsheaf $ \cS'
\subset \cF$ of the same rank as $\cS$ such that $\cS \subset \cS'$
and $\cS \ne \cS'$.
\endclaim

If $\omega$ is a K\"ahler form on $X$ and if $\cS$ is the
$\omega$-maximal destabilizing subsheaf of the $\omega$-unstable sheaf
$\cF$, then $\cS$ is maximal. This is the way we will identify maximal
subsheaves.

\claim 6.12\ Proposition| Let $X$ be a normal compact K\"ahler space
with $a(X) = 0$ and $\cF$ a reflexive coherent sheaf on $X$. Then
$\cF$ admits only finitely many maximal reflexive subsheaves of rank~$1$.
\endclaim

\proof.  Of course we may assume $X$ smooth.  Consider now the maximal
subsheaves $\cS_i \subset \cF$ of rank 1, $i \in I = \bN$. Choose $m
\in \bN$ and $i_1 < \ldots < i_{m}$ such that
$$
\cS' = \cS_{i_1} + \ldots + \cS_{i_{m}} \subset \cF$$
has the
following property: if $j$ is different from the $i_j$, then $\rk \cS'
= \rk (S'+S_j)$.  So things come down to show that there are only
finitely many $j$ such that
$$
\rk (\cS_j \cap \cS') = 1.
$$
In order to prove this, we assume to the contrary that there are infinitely
many $j$ such that $ \rk (\cS_j \cap \cS') = 1$. Then we have
infinitely many subsheaves
$$
\cT_j := \cS_j \cap \cS' \subset \cS'
$$ of rank 1 (use again the finiteness of hypersurfaces in $X.$) 
Now fix $j_0$.  Then by (6.10) there are only finitely many $j$ 
such that $\cT_j \not \subset \cT_{j_0}$.
For all others we have $\cT_j \subset \cT_{j_0}$ and for those we write 
$$
\cT_j = \cS_j - A_j
$$
and
$$
\cT_{j_0} = \cS_{j_0} - A_{j_0}
$$
with effective divisors $A_j$. Since $X$ has only finitely many
irreducible hypersurfaces, we have $A_{j_0} \subset A_j$ for almost
all $j$, hence we obtain $\cS_j \subset \cS_{j_0}$ for almost all $j$,
contradiction to maximality.\qed

\claim 6.13\ Corollary| Let $X$ be an normal compact K\"ahler space
with $a(X) = 0$ and $\cF$ a torsion free sheaf of rank at most $3$.
Then $\cF$ contains only finitely many maximal reflexive subsheaves.
\endclaim

\proof. By 6.12 we have only to deal with the case of subsheaves of
rank $2$. This is done by dualizing and applying 6.12 to $\cF^*$ using
the following trivial remark: if $\cS \subset \cF$ is maximal with
quotient $Q$, then $Q^* \subset \cF^*$ is maximal.\qed
\vskip 20pt

\section{\S7  An abundance theorem for K\"ahler threefolds.}

Here we want to solve (the remaining part of) the abundance problem
for K\"ahler threefolds:

\claim 7.1\ Theorem| Let $X$ be a $\bQ$-Gorenstein K\"ahler threefold
with only terminal singularities such that $K_X$ is nef (a minimal
K\"ahler threefold for short).Then $\kappa (X) \geq 0$.
\endclaim

Of course, more should be true:

\claim 7.2\ Conjecture| Let $X$ be a minimal K\"ahler threefold. Then
$K_X$ is semi-ample, i.e.\ some multiple $mK_X$ is spanned by global
sections.  \endclaim

\claim 7.3\ Remark|{\rm 
{\item {\rm(1)} In case $X$ is projective,
 everything is proved by Miyaoka [Mi87,88] and Kawamata [Ka92].
\item {\rm(2)} In the non-algebraic case, 7.1/7.2 is proved in [Pe00] with 
the important possible exception that $X$ is simple and not Kummer (see 1.4).
In particular in this remaining case we have algebraic dimension $a(X) = 0$
and $\pi_1(X)$ finite. In [DPS00] 7.1 is proved if $K_X$ carries a sufficiently
nice metric, e.g.\ if $K_X$ is hermitian semipositive. 
\item {(3)} In case $X$ is Gorenstein, we have the Riemann-Roch formula
$$ \chi(X,\cO_X) = - {{1} \over {24}} K_X \cdot c_2(X).$$
Therefore the inequality (6.1) -- recall we may assume that $a(X) = 0$ --
$$
K_X \cdot c_2(X) \geq 0 \eqno (*)
$$
implies $\chi(X,\cO_X) \leq 0$ and therefore $h^0(X,K_X) = h^3(X,\cO_X) \ne 0$,
so that at least $\kappa (X) \geq 0$. 
In case $X$ is not Gorenstein, this Riemann-Roch formula is not true;
instead one has some positive correction term [Fl87] which might 
correct the negativity of $-K_X \cdot c_2(X)$ and therefore destroy
the contradiction.}}
\endclaim

\proof.  {\it of Theorem 7.1} \ As noticed in 7.3 we may assume that
$X$ is simple non-Kummer, in particular $q(X) = 0$.  First we reduce
ourselves to the case that $X$ is $\bQ$-factorial by applying
Kawamata's factorialisation $f: \hat X \la X$ as in the proof of 4.4.
Since $f$ is small, we have $K_{\hat X} = f^*(K_X)$, so $K_{\hat X}$
is nef. Hence we can work on $\hat X$ and thus may assume $X$ to be
$\bQ$-factorial from the beginning.  \smallskip
\noindent  We consider 
a desingularization
$$
\pi: \hat X \la X
$$
and compute by Riemann-Roch
$$
\chi(\hat X, \pi^*(mK_X)) = {{m} \over {12}} K_X \cdot c_2(X) +
\chi(X,\cO_X) \eqno (*)
$$
for all $m$ such that $mK_X$ is Cartier.  Assume $\kappa (X) = -
\infty$, so $H^3(\cO_X) = 0$. Since $X$ is not projective, we have
$H^2(\cO_X) \ne 0$. In total we obtain:
$$
\chi(X,\cO_X) \geq 2.
$$
If now $K_X \ne 0$, then by 4.3/5.1, we have $h^2(X,mK_X) \leq 1$,
hence (6.1) and $(*)$ imply $h^0(X,mK_X) \geq 1$, a contradiction.  If
however $K_X = 0$, take a positive integer $m$ such that $mK_X$ is
Cartier. If now $mK_X$ is not a torsion line bundle, we must have
$q(X) > 0$, contradiction.  \qed

\claim 7.4 Remark| {\rm In order to settle the abundance for K\"ahler
threefolds completely, it remains to show that a simple threefold
$X$ with $K_X$ nef and $\kappa (X) = 0$ must be Kummer. In the
following we collect what we know about $X$. We shall assume that
$q(X) = 0$, otherwise we consider the Albanese and are easily done.
Thus we have $\chi(X,\cO_X) \geq 1$.  \smallskip \noindent (1) $K_X
\cdot c_2(X) = 0$ and $1 \leq \chi(X,\cO_X) \leq 2$. \smallskip
\noindent The first part follows easily from equation $(*)$ in the
proof of (7.3) together with 4.1/5.1. Hence
$$
\chi(X,mK_X) = \chi(X,\cO_X) \leqno (*)
$$
for all integers $m$ such that $mK_X$ is Cartier. Then again
4.3/5.1 gives the inequality for $\chi(X,\cO_X)$.\smallskip 

\noindent 
(2) $X$ cannot be Gorenstein. In fact, then the Riemann-Roch formula
$$
24 \chi(X,\cO_X) = -K_X \cdot c_2(X) = 0
$$
gives a contradiction.
\smallskip \noindent
(3) If $ \chi(X,\cO_X) = 2$, then $K_X^2 = 0$. This is a consequence 
of (1) via the vanishing 4.3.
\smallskip \noindent
(4) If $\chi(X,\cO_X) = 1$, then $h^0(X,K_X) = 1$ and $h^2(X,\cO_X) = 1$. }
\endclaim

\section{\S8 Almost algebraic K\"ahler threefolds} 

In this section we show that simple non-Kummer threefolds are very far
from projective threefolds, in a sense which is made precise in the
following definition.

\claim 8.1\ Definition| Let $X$ be a normal K\"ahler variety with only
terminal singularities.  $X$ is almost algebraic if there exists an
algebraic approximation of $X$. This is a proper surjective flat
holomorphic map $\pi: \cX \to \Delta$ from a normal complex space
$\cX$ where $\Delta \subset {\Bbb C}^m$ is the unit disc, where $X
\simeq X_0$, where all complex analytic fibers $X_t = \pi^{-1}(t)$ are
normal K\"ahler spaces with at most terminal singularities such that
there is a sequence $(t_j)$ in $\Delta$ converging to $0$ so that all
$X_j := X_{t_j}$ are projective. 
\endclaim
 
Of course, in case $X$ is smooth, all $X_t$ will be smooth (after
possibly shrinking $\Delta)$.

The following problem is attributed to Kodaira.

\claim 8.2\ Problem|  Is every compact K\"ahler manifold almost algebraic?
\endclaim

From a point of view of algebraic geometry almost algebraic K\"ahler
spaces seem to be the most interesting K\"ahler spaces. Therefore it
is worthwile to notice

\claim 8.3\ Theorem| Let $X$ be a nearly algebraic K\"ahler threefold
with only terminal singularities.  If $X$ is simple and additionally
$K_X$ nef or $X$ smooth, then $X$ is Kummer.
\endclaim

\proof.  Assume that $X$ is not Kummer. Then $\pi_1(X)$ is finite by
[Ca94] as already mentioned.  Let $\pi: \cX \to \Delta$ be an
algebraic approximation of $X$. Let $(t_j)$ be a sequence in $\Delta$
converging to $0$ such that all $X_j = X_{t_j}$ are projective. Notice
first that $\kappa (X_j) \geq 0$ for all $j$. In fact, otherwise $X_j$
would be uniruled for some $j$ and by standard arguments $X_t$ would
be uniruled for all $t$ which is not possible, $X = X_0$ being simple.
\smallskip

\noindent (1) We show that $\kappa (X) = \kappa
(X_0) \geq 0$. Fix a positive integer $m$.
Then by [KM92, 1.6], every $t_j$ admits an open neighborhood $U_j$
such that
$$
h^0(X_t,mK_{X_t}) = h^0(X_j,mK_{X_j}) $$
for all $t \in U_j$. Now choose
$m$ such that $h^0(X_j,mK_{X_j}) > 0$ for some $j$. Then it follows that
$h^0(X_t,mK_{X_t}) = h^0(X_j,mK_{X_j}) =: d > 0$ for all $t$ in an open set 
in $\Delta$.  Let
$$
A := \{t \in \Delta \vert h^0(X_t,mK_{X_t}) \geq d \}.
$$
Then $A$ is an analytic set in $\Delta$ (semi-continuity in the
analytic Zariski topology), and it contains a non-empty open set,
hence $A = \Delta$.  Thus $\kappa (X_0) \geq 0$. Since $X_0$ is
simple, we conclude $\kappa (X_0) = 0$.  \smallskip
\noindent (2) Suppose that $\kappa (X_j) \geq 1$
for some $j$. Then fix $m$ such that $h^0(X_j,mK_{X_j}) \geq 2$. Repeating
the same arguments as in (1), we conclude $h^0(X,mK_X) \geq 2$,
contradicting $X$ being simple. So $\kappa (X_j) = 0$ for all $j$. 
\smallskip
\noindent
(3) Here we will show that $X_j$ is Kummer for all $j$. 
Let
$X'_j$ be a minimal model of $X_j$. Observe that
$$
h^2(X_j,\cO_{X_j}) = h^2(X,\cO_X) > 0;
$$
in fact, $H^1(X,\cO_X) = 0$, hence $H^1(X_j,\cO_{X_j}) = 0 $ for large
$j$. Moreover $h^0(X_t,K_{X_t})$ is constant by [KM92], as shown above.
Therefore the equality follows by Serre duality and the constancy of
$\chi(X_t,\cO_{X_t})$.  Hence we also have $h^2(X'_j,\cO_{X'_j}) > 0$.
\smallskip

\noindent Since $K_{X'_j} \equiv 0$, there exists a finite cover, the
so-called canonical cover, $h: \tilde X_j \to X'_j$, \'etale in
codimension $2$, such that $K_{\tilde X_j} = \cO_{\tilde X_j}$. In
particular $\tilde X_j$ is Gorenstein and Riemann-Roch yields
$$
\chi(\tilde X,\cO_{\tilde X_j}) = 0.
$$
Since $h^2(\cO_{\tilde X_j}) > 0$, we must have $q(\tilde X_j) >
0$. .  Let $\alpha_j: \tilde X_j \la A = A_j$ be the Albanese map. By
[Ka85] there exists a finite \'etale cover $B \la A$ such that
$$
\hat X_j := \tilde X_j \times_A B \simeq F \times B.
$$
In particular $\hat X_j$ and $\tilde X_j$ are smooth because of the
isolatedness of singularities.  We conclude that $X_j'$ is Kummer
unless $F$ is a K3-surface. To exclude that case, consider the image
$F' \subset \tilde X_j'$ of a general $F \times \{b\}$. Then $F'$ is
K3 or Enriques and does not meet the singularities of $X_j'$. Moreover
the normal bundle $N_{F'}$ is numerically trivial. Since $F'$ moves,
it is actually trivial. Now consider the strict transform in $X_j$,
again called $F'$. Then $F'$ has the same normal bundle in $X_j$, so
that $N_{F'/\cX} = \cO_{\cX}^2$. Since $H^1(N) = 0$, the deformations
of $F'$ cover every $X_t$ contradicting the simplicity of $X_0$. So
$F$ cannot be K3 and $X_j$ is Kummer.  \smallskip

\noindent (4) Suppose $K_{X_0} $ nef.  Fix a positive number $m$ such that 
$$
mK_{\cX} = \cO_{\cX}(D)
$$
with some effective divisor $D$.  We may assume that $D$ does not
contain any fiber of $\pi;$ denote $D_t = X_t \cap D$.  We want to
argue that $K_{X_t}$ must be nef, therefore $K_{X_j} = 0$ so that $D_j
= 0$ and $D = 0$ in total. So we will obtain $K_X = 0$. To see that
$X$ is Kummer, consider the canonical cover of $X$ and argue as in the
proof of (7.1). To prove nefness, we apply [KM92] to deduce that the
sequence $\varphi_j: X_j \to X_j' = A_j/G$ appears in family
$\cX_{U_j} \la \cX'_{U_j} $ over a small neighborhood $U_j$ of $t_j$.
In particular some multiple $N^{*\mu}_{D_t}$ has many sections for $t
\in U_j$ on an at least 1-dimensional family of curves. By
semicontinuity, also $N^{*\mu}_{D_0}$ has many sections on such a
family, contradicting the nefness of $D_0$.  Alternatively, $K_{ \cX}$
is negative on a family of rational curves over $U_j$, which converges
to a family of rational curves in $X_0$ and therefore forces $K_{X_0}$
to be non-nef.  \smallskip

\noindent (5) Now suppose $X_0$ smooth, i.e. $\pi$ is smooth after 
shrinking $\Delta$.  Take a sequence of blow-ups of smooth
subvarieties of $\cX$ such that the preimage of ${\rm red} D$ has
normal crossings. After shrinking $\Delta$ we may assume that the only
points and compact curves blown up lying over $X_0$ so that all fibers
over $\Delta \setminus 0$ are smooth.  Then take a covering $h: \tilde
\cX \la \cX$ such that $K_{\tilde \cX} = \cO(\tilde D)$ with $\tilde
\cX$ smooth. This is possible e.g. by applying [Ka81]. Then
$\cX_{t_j}$ is Kummer and admits a 3-form and therefore must be
bimeromorphically a torus (if $A/G$ admits a 3-form, then it is a
torus covered by $A$. This is a consequence of the simplicity of $A$
and the fact that $G$ acts without fix points).  Hence every $\tilde
\cX_t$, $t \ne 0$, has 3 holomorphic 1-forms which are independent at
the general point and therefore every $X_t$, $t \ne 0$, is Kummer. In
order to show that $X_0$ is Kummer, consider the central fiber $\tilde
\cX_0$ which contains the preimage of the strict transform $X'_0$ of
$X_0$. More precisely, we have
$$
\tilde \cX_0 = X'_0 + \sum a_i E_i
$$
where the $E_i$ are smooth threefolds contracted to points or
curves.  By semi-continuity, $h^2(\cO_{\tilde \cX_0}) \geq 3$. Now we
check easily that
$$
H^2(\tilde\cX_0,\cO_{\tilde \cX_0}) = H^2(X'_0,\cO_{X'_0}),
$$
hence $X'_0$ carries three 2-forms coming from $\tilde \cX$. But
then it is clear that also some of the holomorphic 1-forms on $\tilde
\cX$ give non-zero 1-forms on $X'_0$, since the 2-forms are wegdge
products of the 1-forms. Hence $X'_0$ is Kummer and so does $X_0$.
\qed

\bigskip

\noindent 
If problem 8.2 has a positive answer in dimension 3, Theorem 8.3 
excludes the existence of simple non-Kummer threefolds.

\section{References}

\bigskip

{\eightpoint
  
\bibitem [BS95]&Beltrametti, M.;Sommese, A.J.:& The adjunction theory of
  complex projective varieties.& de Gruyter Exp. in Math. {\bf 16}
  (1995)&
  
\bibitem [Ca94]&Campana, F.:& Remarques sur le rev\^etement universel des
  vari\'et\'es K\"ahl\'eriennes compacts;& Bull.\ Soc.\ Math.\ France
  {\bf 122} (1994) 255--284&

\bibitem [CP01]&Campana, F.; Peternell, Th.:& The Kodaira dimension of Kummer 
 threefolds;& Bull.\ Soc.\ Math.\ France {\bf 129} (2001) 357--359&

\bibitem [De92]&Demailly, J.-P.:& Regularization of closed
  positive currents and intersection theory;& J.\ Alg.\ Geom.\ {\bf 1} (1992),
  361--409&
  
\bibitem [De93a]&Demailly, J.-P.:& Monge-Amp\`ere operators, Lelong
  numbers and intersection theory;& Complex Analysis and Geometry,
  Univ.\ Series in Math., edited by V.~Ancona and A.~Silva, Plenum
  Press, New-York (1993)&

\bibitem [De93b]&Demailly, J.-P.:& A numerical criterion for very ample line 
  bundles;& J.~Differential Geom.\ {\bf 37} (1993) 323--374&

\bibitem [DPS94]&Demailly, J.-P.; Peternell, Th.; Schneider, M.:&
  Compact complex manifolds with numerically effective tangent
  bundles;& J.\ Alg.\ Geom.\ {\bf 3} (1994) 295--345&

\bibitem [DPS00]&Demailly, J.-P.; Peternell,Th.; Schneider,M.;& 
  Pseudo-effective line bundles on compact K\"ahler manifolds;& 
  Intern.\ J.\ Math. {\bf 6} (2001) 689--741&
  
\bibitem [Ei95]&Eisenbud, D.: & Commutative algebra;& Graduate Texts in
  Math.\ {\bf 150}, Springer (1995)&
  
\bibitem [En87]&Enoki, I.:& Stability and negativity for tangent
  bundles of minimal K\"ahler spaces;& Lecture Notes in Math. {\bf
  1339} (1987) 118--127&

\bibitem [Fl87]&Fletcher, A.R.:& Contributions to Riemann-Roch on
  projective $3$-folds with only canonical singularities and
  applications;& Proc.\ Symp.\ Pure Math.\ {\bf 46} (1987) 221--231&

\bibitem [Gr62]&Grauert, H.:& \"Uber Modifikationen und exzeptionelle
  analytische Mengen;& Math.\ Ann.\ {\bf 146} (1962) 331--368&

\bibitem [GR70]&Grauert, H.; Riemenschneider,O.:&
  Verschwin\-dungs\-s\"atze f\"ur analy\-tische                          
  Koho\-mo\-lo\-gie\-grup\-pen auf komplexen R\"aumen;& Invent.\ Math.\
  {\bf 11} (1970) 263--292&

\bibitem [Ka81] &Kawamata, Y.:& Characterization of abelian varieties.& 
  Comp. \ math. \ {\bf 43} (1981) 275-276&

\bibitem [Ka85] &Kawamata, Y.:& Minimal models and the Kodaira dimension of
  algebraic fiber spaces;& J. \ reine u. angew. \ Math. \ {\bf 363} 
  (1985) 1-46& 

\bibitem [Ka88]&Kawamata, Y.:& Crepant blowing ups of threedimensional
  canonical singularities and applications to degenerations of
  surfaces;& Ann.\ Math.\ {\bf 119} (1988) 93--163&

\bibitem [KMM87]&Kawamata, Y.; Matsuki, K.; Matsuda, K.:& Introduction
  to the minimal model program;& Adv.\ Stud.\ Pure Math.\ {\bf 10}
  (1987) 283--360&

\bibitem [KM92]&Koll\'ar, J.; Mori, S.:& Classification of
  three dimensional flips;& Journal of the AMS {\bf 5} (1992), 533--703&

\bibitem [Ko87]&Kobayashi, S.:& Differential geometry of complex vector
  bundles;& Princeton Univ.\ Press (1987)&
  
\bibitem[Ko92]&Kolla\'r, J. et al.:& Flips and abundance for algebraic
  $3$-folds;& Ast\'erisque {\bf 211}, Soc.\ Math.\ France (1992)&
  
\bibitem [Mi87]&Miyaoka, Y.:& The Chern classes and Kodaira dimension
  of a minimal variety;& Adv.\ Stud.\ Pure Math.\ {\bf 10} (1987) 
   449--476&

\bibitem [Mi88]&Miyaoka, Y.:& On the Kodaira dimension of minimal
  threefolds;& Math.\ Ann.\ {\bf 281} (1988) 325--332 &

\bibitem [Mi88a]&Miyaoka, Y.:&  Abundance conjecture for 3-folds: case 
  $\nu = 1$;& Comp.\ Math.\ {\bf 68} (1988) 203--220& 

\bibitem [Pa98]&Paun, M.:& Sur l'effectivit\'e num\'erique des
  images inverses de fibr\'es en droites;& Math.\ Ann.\ {\bf 310}
  (1998) 411--421&
  
\bibitem [Pe01]&Peternell, Th.:& Towards a Mori theory on compact
  K\"ahler $3$-folds, III;& Bull.\ Soc.\
  Math.\ France {\bf 129} (2001) 339--356&
  
\bibitem [Re87]&Reid, M.:& Young person's guide to canonical
  singularities;& Proc. Symp. Pure Math. {\bf 46}, 345-414 (1987)&

\bibitem [Siu74]&Siu, Y.T.:& Analyticity of sets associated to Lelong 
  numbers and the extension of closed positive currents;& Invent.\ Math.\
  {\bf 27} (1974) 53--156&

\bibitem [Va84]&Varouchas, J.:& Stabilit\'e de la classe des vari\'et\'es 
  K\"ahl\'eriennes par certain morphismes \  propres;& 
  Invent.\ Math.\ {\bf 77}, (1984)  117--127&

}

\end